 \documentclass[10pt,a4paper]{article}

\usepackage{graphicx, color, epstopdf, amsmath,subfigure, amsfonts, hyperref} 

\catcode`\@=11
\@addtoreset{equation}{section}

\catcode`\@=12

\newtheorem{Theorem}{Theorem}[section]
\newtheorem{Proposition}{Proposition}[section]
\newtheorem{Lemma}{Lemma}[section]
\newtheorem{Corollary}{Corollary}[section]
\newtheorem{Remark}{Remark}[section]

\def\hb{\hfill\break}
\def\fa{\forall}
\def\Proof{\noindent{{\bf Proof. }}}
\def\square{\vbox{ \hrule height .4pt
            \hbox{\vrule width .4pt height 7pt \kern 7pt
			\vrule width .4pt}\hrule height .4pt }}
\def\QED{\hfill {$\square$}\goodbreak \medskip}

\def\norm#1{\| #1\|}

\def\f{{V}}
\def\bs{{\bar\sigma}}
\def\bt{{\bar\tau}}
\def\ff{{\varphi}}

\def\N{{\mathbb N}}

\def\Z{{\mathbb Z}}
\def\R{{\mathbb R}}

\def\aub{{\alpha_{u,b}}}
\def\oub{{\omega_{u,b}}}
\def\g{{\gamma}}
\def\e{\varepsilon}

\def\ub{{\beta}}
\def\bi{{\beta}}

\def\s{{\sigma}}

\def\VV{{\cal V}}

\def\MM{{\cal H}}

\def\epsilon{{\varepsilon}}

\def\XX{{\cal X}}
\def\X{{\cal X}}

\def\dist{{\mathrm{dist}}}

\title{{ An energy constrained method for the existence of layered type solutions of NLS equations.}}
\author{Francesca Alessio$^1$
and Piero Montecchiari$^2$}
\date{ }

\begin{document}

\smallskip

\maketitle
\centerline{\small{
 Dipartimento di Ingegneria Industriale e Scienze
Matematiche, 
} 
}
\centerline{
\small{Universit\`a Politecnica delle Marche, Via Brecce Bianche -- I 60131 Ancona,
}}
\centerline{
\small{
e-mail:
$^{1}${\tt f.g.alessio@univpm.it}, $^{2}${\tt p.montecchiari@univpm.it}.
}
}

\bigskip\bigskip
\bigskip\bigskip

\noindent{\small{\bf Abstract.}
We study the existence of
positive solutions on $\R^{N+1}$ to semilinear elliptic equation
$
-\Delta u+u=f(u)
$
where $N\geq 1$ and $f$ is modeled on the power case $f(u)=|u|^{p-1}u$.
Denoting with $c$ the mountain pass level of $\f(u)=\tfrac 12\|u\|^{2}_{H^{1}(\R^{N})}-\int_{\R^{N}}F(u)\, dx$, $u\in H^{1}(\R^{N})$ ($F(s)=\int_{0}^{s}f(t)\, dt$), we show, via a new energy constrained variational argument, that for any $b\in [0,c)$ there exists a positive bounded solution $v_{b}\in C^{2}(\R^{N+1})$ such that $E_{v_{b}}(y)=\tfrac 12\|\partial_{y}v_{b}(\cdot,y)\|^{2}_{L^{2}(\R^{N})}-V(v_{b}(\cdot,y))=-b$ and $v(x,y)\to 0$ as $|x|\to+\infty$ uniformly with respect to $y\in\R$. We also characterize the monotonicity, symmetry and periodicity properties of $v_{b}$.}

\bigskip\bigskip


\bigskip


\vskip3truecm
\noindent{
{\scriptsize $^{1,2}$ Partially supported by the PRIN2009 grant "Critical Point Theory and Perturbative Methods for Nonlinear Differential Equations"}}

\vfill\vfill\eject

%

\section{Introduction}

In this paper we study the existence of positive solutions on $\R^{N+1}$ to semilinear elliptic equations
$$
-\Delta u+u=f(u)\eqno(E)
$$
where $N\geq 1$ and $f$ is a nonlinearity which can be thought modeled on the power case $f(u)=|u|^{p-1}u$ with $p$ subcritical and greater than 1. Equations of this kind are used in various fields of Physics such as, for example, plasma or laser self-focusing models (see \cite{[SS]} and the references therein). They arise in particular in the study of standing waves (stationary states) solutions of the corresponding nonlinear Schr\"{o}dinger type equations.
 
Starting with the work by W. A. Strauss, \cite{[S]}, the problem of finding and characterizing positive solutions $v\in H^{1}(\R^{N+1})$ of (E) has been widely studied. We refer to the paper by H. Berestycki and P.L. Lions \cite{[BL]} (in the case $N\geq 2$, see \cite{[BGK]}  for $N=1$) where nearly optimal existence results regarding least energy solutions for (E) are obtained. Their mountain pass characterization, and so information about their Morse index, is given by L. Jeanjean and K. Tanaka in \cite{[JJT]}. In the pure power case, uniqueness and non degeneracy properties of solutions of (E) in $H^{1}(\R^{N+1})$ was derived by M.K. Kwong in \cite{[K]}. Regarding the uniqueness problem for more general nonlinearity $f$, we refer to the paper by J. Serrin and M. Tang, \cite{[ST]}, and to the references therein.

A new kind of entire solutions of (E) has been introduced  by N. Dancer in \cite{[D]}. Denoting $(x,y)\in\R^{N}\times\R$ a point in $\R^{N+1}$, we note that a ground state solution $u_{0}(x)$ of (E) in $\R^{N}$ can be  thought as a solution of (E) on $\R^{N+1}$, which  is constant with respect to the $y$ variable. In the pure power case (or anyhow assuming the nondegeneracy of the ground state solution)  Dancer proved, by using bifurcation and continuation arguments, the existence of a continuous branch of entire positive solutions of (E) in $\R^{N+1}$ bifurcating from the {\sl cylindric type} solution $u_{0}$. These solutions are periodic in the variable $y$ and decay to zero as $|x|\to+\infty$. Different periodic Dancer's solutions (suitably rotated) were then used in the pure power case as prescribed asymptotes in the constructions of {\sl multiple ends} solutions of (E) by A. Malchiodi in \cite{[Ma]} and by M. del Pino, M. Kowalczyk, F, Pacard and J. Wei in \cite{[DpKPW]}. 

Related to the above papers is  the one by C. Gui, A. Malchiodi and H. Xu, \cite{[GMH]}, where qualitative properties (such as radial symmetry with respect to the variable $x$ and eveness with respect to $y$)  of positive solutions $v(x,y)$ of $(E)$ which decay to zero as $|x|\to+\infty$ (uniformly w.r.t. $y$) are established. Their study is based on moving plane techniques together with the use of some Hamiltonian identities  which are connected with the Lagrangian structure of that kind of problem. 

To describe the Hamiltonian identities which are used in \cite{[GMH]} and to introduce precisely the problem studied in the present paper, note that prescribing the decay properties of a solution $v$ only with respect to the variable $x\in\R^{N}$, naturally gives to the variable $y$ the role of an evolution
variable. In this respect, as usual in the evolution problems, all the solutions $v$ of (E) described above belong to the
space $X=L^{2}_{loc}(\R, H^{1}(\R^{N}))\cap H^{1}_{loc}(\R, L^{2}(\R^{N}))$ and verify (at least in a weak sense) the evolution equation
\begin{equation}\label{eq:eq}
\partial_{y}^{2}v(\cdot,y)= V'(v(\cdot,y)),\quad y\in\R,
\end{equation}
where $V'$ is the gradient in $H^{1}(\R^{N})$ of the Euler functional relative to equation (E) on $\R^{N}$,
$$
V(u)=\int_{\R^{N}}\tfrac 12|\nabla u|^{2}+\tfrac 12|u|^{2}-F(u)\, dx,\quad u\in H^{1}(\R^{N}),
$$
where $F(s)=\int_{0}^{s}f(t)\, dt$.
We will refer to this kind of solutions  as {\sl layered solutions}  of (E).   

Noting that  equation (\ref{eq:eq}) has Lagrangian structure, one can think to the variable $y$ as a {\sl time} variable and to the functional $U=-V$ as the {\sl energy potential} of the infinite dimensional  dynamical system. Every layered solution $v$ defines a trajectory $y\in\R\to v(\cdot,y)\in H^{1}(\R^{N})$, solution to (\ref{eq:eq}).
 In this connection, any $u\in H^{1}(\R^{N})$ which solves (E) is an equilibrium of (\ref{eq:eq}) and the solutions found by Dancer are periodic orbits of the system. Since the system is autonomous, if $v$ is a  layered solution to (E) then the {\sl Energy} function
$$y\to E_{v}(y)=\frac{1}{2}\|\partial_{y} v(\cdot,y)\|_{L^{2}(\R^{N})}^{2}-V(v(\cdot,y))$$
is constant (a formal proof of this Hamiltonian identity for a general class of elliptic equations can be found in \cite{[BF]} and \cite{[G]}, see also \cite{[AM1]} for the case of Allen Cahn equations).   \bigskip
 
 In the present paper, in analogy with the study already done for Allen Cahn type equation in \cite{[AM1]}, \cite{[AM2]}, \cite{[AM3]} (see also  \cite{[A]} for Allen Cahn system of equations),
 we study the problem of finding layered solution of (E) with prescribed energy. In particular we study the problem of looking for {\sl connecting orbit} solutions with prescribed energy.\smallskip
 
 To be more detailed, we precise our assumption on the non linearity $f$. We assume that 
\begin{description}
\item[$(f1)$] $f\in C^1(\R)$,
\item[$(f2)$] there exists $C>0$ and $p\in(1,1+\tfrac{4}{N})$ such that
$|f(t)|\le C(1+|t|^{p})$ for any $t\in\R$,
\item[$(f3)$] there exists $\mu>2$ such that $0<\mu F(t)\le f(t)t$
for any $t\ne 0$, where $F(t)=\int_0^t f(s)\,ds$,
\item[$(f4)$] $f(t)t<f'(t)t^{2}$ for any $t\ne 0$.
\end{description}
As it is well known, $(f1)$--$(f4)$ are more than sufficient to guaranty that $V\in C^{1}(H^{1}(\R^{N}))$ and that it satisfies the geometrical assumptions of the Mountain pass Theorem. Setting $c=\inf_{\gamma\in\Gamma}\sup_{t\in[0,1]}V(\gamma(t))$, where $\Gamma=\{\gamma\in C([0,1],H^{1}(\R^{N}))\,|\, \gamma(0)=0,\, V(\gamma(1))<0\}$, we have that $c>0$ is an asymptotical critical level for $V$. Concentration compactness arguments allow to prove that $c$ is actually the lowest positive critical level of $V$. Then, the definition of the Mountain pass level implies that given any
$b\in [0,c)$ the sublevel $\{V\leq b\}$ is the union of two disjoint path connected sets ${\cal V}^{b}_{-}$ and ${\cal V}^{b}_{+}$, where we denote with ${\cal V}^{b}_{-}$ the one which contains $0$. The main result of the present paper establishes that given any $b\in [0,c)$ there exists a layered solution $v$ of (E) with $E_{v}=-b$ and which connects the set ${\cal V}^{b}_{-}$ and
${\cal V}^{b}_{+}$, in the sense that 
$\liminf_{y\to\pm\infty}\dist_{L^{2}(\R^{N})}(v(\cdot,y),{\cal V}^{b}_{\pm})=0$.
Precisely we prove that

\begin{Theorem}\label{T:main} If $F$ satisfies $(f1)-(f4)$ then for any $b\in [0,c)$ the equation (E) has a solution $v_{b}\in C^{2}(\R^{N+1})$ with {\sl energy} $E_{v_{b}}=-b$ and such that
\begin{description}
\item{i)} $v_{b}>0$ on $\R^{n+1}$, 
\item{ii)} $v_{b}(x,y)=v_{b}(|x|,y)\to 0$ as $|x|\to+\infty$, uniformly w.r.t. $y\in\R$,
\item{iii)} $\partial_{r}v_{b}(x,y)<0$  for $r=|x|>0$ and $y\in\R$.
\end{description}
 Moreover, if $b>0$, 
\begin{description}
\item{iv)} there exists $T_{b}>0$ such that $v_{b}$ is periodic of period $2T_{b}$ in the variable $y$ and symmetric with respect to $y=0$ and $y=T_{b}$.
\item{v)} $\partial_{y}v_{b}(x,y)>0$ on $\R^{N}\times (0,T_{b})$, $v_{b}(\cdot,0)\in\VV^{b}_{-}$, $v_{b}(\cdot,T_{b})\in\VV^{b}_{+}$.
\end{description}
Finally, if $b=0$, 
\begin{description}
\item{v)}  $v_{0}\in H^{1}(\R^{N+1})$ is radially symmetric and $\partial_{r}v_{0}<0$ for $r=|(x,y)|>0$,
\item{vi)}  $v_{0}(\cdot,0)\in\VV^{0}_{+}$ and $v_{0}$ is a mountain pass point of  the Euler funcional relative to (E) on $H^{1}(\R^{N+1})$.
\end{description}
\end{Theorem}
Theorem \ref{T:main} gives the existence for any $b\in [0,c)$ of a positive layered solution $v_{b}$ to (E) with energy $-b$ which is radially symmetric and decaying to $0$ as $|x|\to+\infty$ uniformly with respect to $y\in\R$.  When $b>0$ the solution $v_{b}$ is a {\sl periodic solution} of period $2T_{b}$ which is symmetric with respect to $y=0$ and $y=T_{b}$. It can be thought as a trajectory which oscillates back and forth along a simple curve connecting the two turning points $v_{b}(\cdot,0)\in{\cal V}^{b}_{-}$ and $v_{b}(\cdot,T_{b})\in{\cal V}^{b}_{+}$. These solutions, which we call  {\sl brake orbit type solutions}, are clearly analogous to the Dancer solutions. When $b=0$ the solution $v_{0}$ defines a trajectory which emanates from $0\in H^{1}(\R^{N})$ as $y\to-\infty$, reaches the point $v(\cdot,0)\in {\cal V}^{0}_{+}$ and goes back symmetrically to $0$ for $y>0$. It can been thought as a {\sl homoclinic solution} to $0\in H^{1}(\R^{N})$ and it is in fact the mountain pass point of  the Euler funcional relative to (E) on $H^{1}(\R^{N+1})$. Finally we can think at the mountain pass point of $V$ in $H^{1}(\R^{N})$ as an equilibrium of (\ref{eq:eq}) at energy $-c$. The Energy diagram here below wants to summarize these considerations.

\bigskip

\noindent
\def\svgwidth{0.8\columnwidth}

\begingroup
  \makeatletter
  \providecommand\color[2][]{%
    \errmessage{(Inkscape) Color is used for the text in Inkscape, but the package 'color.sty' is not loaded}
    \renewcommand\color[2][]{}%
  }
  \providecommand\transparent[1]{%
    \errmessage{(Inkscape) Transparency is used (non-zero) for the text in Inkscape, but the package 'transparent.sty' is not loaded}
    \renewcommand\transparent[1]{}%
  }
  \providecommand\rotatebox[2]{#2}
  \ifx\svgwidth\undefined
    \setlength{\unitlength}{559.90865568pt}
  \else
    \setlength{\unitlength}{\svgwidth}
  \fi
  \global\let\svgwidth\undefined
  \makeatother
  \begin{picture}(1,0.24577694)%
    \put(0,0){\includegraphics[width=\unitlength]{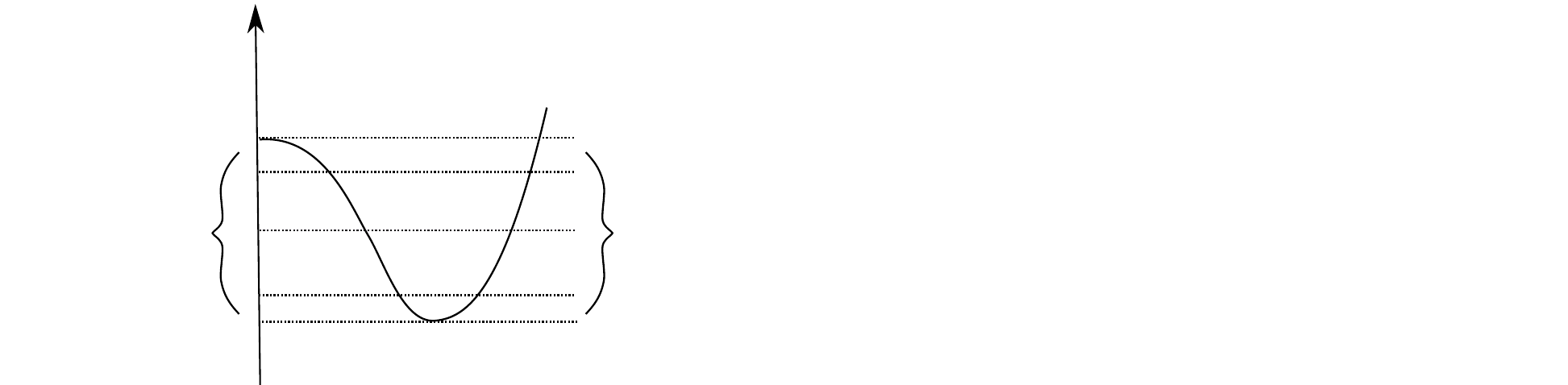}}%
    \put(-0.03218624,0.08867883){\color[rgb]{0,0,0}\makebox(0,0)[lb]{\smash{$-c<E<0$}}}%
    \put(0.404775,0.09528141){\color[rgb]{0,0,0}\makebox(0,0)[lb]{\smash{$\hbox{\small Dancer Solutions}\equiv\hbox{\small Brake orbit type periodic solutions}$}}}%
    \put(0.40503274,0.0339561){\color[rgb]{0,0,0}\makebox(0,0)[lb]{\smash{$\hbox{\small $N$ Mountain Pass point}\equiv\hbox{\small Equilibrium.}$}}}%
    \put(0.40578691,0.15698363){\color[rgb]{0,0,0}\makebox(0,0)[lb]{\smash{$\hbox{\small $N+1$ Mountain Pass point}\equiv\hbox{\small Homoclinic orbit type  solution}$}}}%
    \put(0.17477676,0.2349194){\color[rgb]{0,0,0}\makebox(0,0)[lb]{\smash{$-V$}}}%
    \put(0.05137252,0.15761846){\color[rgb]{0,0,0}\makebox(0,0)[lb]{\smash{$E=0$}}}%
    \put(0.03116758,0.0336984){\color[rgb]{0,0,0}\makebox(0,0)[lb]{\smash{$E=-c$}}}%
  \end{picture}%
\endgroup

\bigskip

To prove Theorem \ref{T:main} we make use of variational methods and we apply an Energy constrained variational argument already introduced and used in \cite{[AM1]}, \cite{[AM2]} and \cite{[AM3]}.
Given  $b\in [0,c)$, we look for minima of the renormalized functional
$$\varphi(v)=\int_{\R} \tfrac{1}{ 2}\|\partial_{y} v
	(\cdot,y)\|^{2}_{L^{2}(\R^{N})}+(V(v(\cdot,y)) -b)\, dy
$$
on the space of function $v\in X$ which are radially symmetric with respect to $x\in\R^{N}$, monotone decreasing with respect to $|x|$ and which verify 
\begin{equation}\label{eq:constraint} \liminf_{y\to
    \pm\infty}\dist_{L^{2}(\R^{N})}(v(\cdot,y),{\cal V}^{b}_{\pm})= 0    \hbox{ and }
     \inf_{y\in\R}V(v(\cdot,y))\geq b.
\end{equation}
Thanks to the constraint $\inf_{y\in\R}V(v(\cdot,y))\geq b$, the functional $\varphi$ is well defined
on this class of functions. Moreover, its minimizing sequences admits limit points $\bar v\in X$ (a priori not verifying (\ref{eq:constraint})) with respect to the weak topology of $H^{1}_{loc}(\R^{2})$. \smallskip

Defining $\bs=\sup\{y\in \R\, /\, v(\cdot,y)\in {\cal V}^{b}_{-}\}$ and $\bt=\inf\{y>\bs\, /\, v(\cdot,y)\in {\cal V}^{b}_{-}\}$, we can prove that  
 $-\infty\leq\bs<\bt<+\infty$ (indeed $\bs>-\infty$ when $b>0$) and
$\lim_{y\to\bs^{+}}\dist(\bar v(\cdot,y),{\cal V}^{b}_{-})=0$, $\bar v(\cdot,\bt)\in{\cal V}^{b}_{+}$ and $V(v(\cdot,y))>0$ for any $y\in (\bs,\bt)$.
Then, the minimality properties of $\bar v$ allow us to prove that 
 $\bar v$ solves in a classical sense the equation (E) on $\R^{N}\times(\bs,\bt)$ and $E_{\bar v}(y)=-b$ for any $y\in (\bs,\bt)$. This will imply that $\bar v$ satisfies the boundary conditions $\lim_{y\to\bs^{+}}\partial_{y}\bar v(\cdot,y)=\lim_{y\to\bt^{-}}\partial_{y}\bar v(\cdot,y)=0$ in $L^{2}$ and
the entire solution $v_{b}$ is  recovered from $\bar v$ 
by translations,  reflections and, eventually, periodic continuations.\bigskip

The variational approach that we used is similar to the one already applyed in the study of the Allen Cahn type equation in \cite{[AM1]}, \cite{[AM2]},\cite{[AM3]}, but the present case is much more complicated due to some natural lack of compactness and weak semicontinuity of the problem. This mainly depends on the competition  between the term $\|u\|_{H^{1}(\R^{N})}^{2}$  and $\int_{\R^{N}}F(u)$ which enters in the definition of the potential functional $\f(u)$ with different sign.  This explain why we assume in $(f2)$ that
$p<1+4/N$. The exponent $p=1+4/N$ is in fact critical  with respect to the existence  for the minimum problem $\inf\{\f(u)\,|\, u\in H^{1}(\R^{N}),\, \|u\|_{L^{2}(\R^{N})}=1\}$ and, related to that, with respect to
the property of orbital stability of the solutions of (E) in $H^{1}(\R^{N})$. We recall that  the ground state
solutions of (E) in $H^{1}(\R^{N})$ are stable when $1<p<1+ 4/N$ (see \cite{[CL]}) while the solutions of (E) in $H^{1}(\R^{N})$ are unstable  
when $1 + 4/N\leq p$ (see \cite{[BC]} for the case $1 + 4/N< p$ and \cite{[W]} for $p=1+4/N$). Another (related) criticality of the exponent $p=1+4/N$ is the fact that the sets ${\cal V}^{b}_{\pm}$ have positive $L^{2}(\R^{N})$ distance if and only if $p<1+4/N$ (one can simply verify it by using dilations in the pure power case). Even if we think that this assumption is only technical and can be probably overcome, here we begin to study the problem in this more compact setting.
\bigskip

The paper is organized as follows. In section 2 we recall some properties of the functional $V$ studying in particular the structure of the sublevel set ${\cal V}^{b}_{\pm}$. The study of the functional $\varphi$ and the use of the energy constraint variational principle described above is contained in section 3.
\bigskip

{\bf Acknowledgments.} We wish to thanks Andrea Malchiodi, Margherita Nolasco and Vittorio Coti Zelati for the useful comments and discussions.
\bigskip

\begin{Remark}\label{R:conseguenze}{\rm Since we look for positive solution of (E) it is not restrictive to assume, and we will do it along the paper, that $f$ is an odd function
\begin{description}
\item[$(f5)$] $f(t)=-f(-t)$ for any $t>0$.
\end{description}
Moreover, we list also some plain consequences of $(f1)$-$(f4)$.
\begin{itemize}
\item[i)] By $(f1)$ and $(f3)$ it is 
straightforward to verify that 
$f(0)=f'(0)=0$ and so
$f(t)=o(t)$ as $t\to 0$.
\item[ii)] By 
(i)
and $(f2)$ we have
\begin{equation}\label{eq:originf}
\fa\e >0,\ \exists A_\e>0\ \hbox{such that}\ |f(t)|\leq \e |t|+ A_\e |t|^p,\
\fa\, t\in\R,
\end{equation}
from which we also derive 
\begin{equation}\label{eq:originF}
\fa\e >0,\ \exists A_\e>0\ \hbox{such that}\ |F(t)|\leq \tfrac{\e}2 |t|^{2}+ \tfrac{A_\e}{p+1} |t|^{p+1},\
\fa\, t\in\R.\end{equation}
\item[iii)] By $(f3)$, if $t\neq 0$ and $s>0$, we have $\tfrac d{ds}F(st)=\tfrac1s f(st)st>\tfrac{\mu}s F(st)$. Hence, 
\begin{equation}\label{eq:superq}F(st)>F(t)s^{\mu}\hbox{ whenever }t\not=0\hbox{ and }s>1.\end{equation}
\item[iv)] By $(f4)$, one plainly verify that, for any $t\neq 0$, 
\begin{equation}\label{eq:convessita}\hbox{the function }s\mapsto\tfrac1s f(st)t\hbox{ is strictly increasing for }s>0.\end{equation}
\end{itemize}}\end{Remark}
For the sake of brevity in the notation,  along the paper we denote $\|u\|\equiv\|u\|_{H^{1}(\R^{N})}$, $\|u\|_{p}=\|u\|_{L^{p}(\R^{N})}$ and $\langle u,v\rangle= \langle u,v\rangle_{H^{1}(\R^{n})}$, $\langle u,v\rangle_{2}= \langle u,v\rangle_{L^{2}(\R^{n})}$ for $n=N$ or $n=N+1$ . Moreover $\dist(A,B)\equiv\dist_{L^{2}(\R^{N})}(A,B)=\inf_{v\in A,\, w\in B}\|v-w\|_{2}$  and
$\dist(u,B)\equiv\inf_{v\in B}\|u-v\|_{2}$ for $A, B\subset L^{2}(\R^{N})$, $u\in L^{2}(\R^{N})$. Given $y\in {\R^N}$ we set $B_{r}(y)\equiv\{x\in\R^{N}\,/\, |x|<r\}$ and $B_{r}\equiv B_{r}(0)$.\\

\section{The Potential functional}
In this chapter, we study some properties of the functional $\f: H^{1}(\R^{N})\to\R$ defined by
\begin{equation}\label{eq:funzionale}
\f(u)=\tfrac12\|u\|^2-\int_{{\R^N}} F(u(x))\,dx.
\end{equation}

\subsection{ The Mountain Pass structure.}
Here we list some classical properties of $\f$, in particular the ones regarding its mountain pass behaviour.\\
First of all we recall that $\f$ is regular on $H^{1}(\R^{N})$ (see e.g \cite{[AP]} and \cite{[M]}). 
\begin{Lemma}\label{L:regV}
 $\f\in {\mathcal C}^2(H^1({\R^N}))$ with
$
\f'(u)h=\int_{\R^{N}} \nabla u\,\nabla h+u\,h- f(u)h\,dx
$ and $\f''(u)h\cdot h=\int_{\R}|\nabla h|^{2}+|h|^{2}-f'(u)h^{2}\, dx$ for all $h\in H^{1}(\R^{N})$.
\end{Lemma}

\noindent Moreover the functional $V$ satisfies the {\sl geometrical hypotheses} of the Mountain Pass Theorem.
Indeed, since $p+1<2^{*}_{N}$, by the Sobolev Immersion Theorem  and Remark \ref{R:conseguenze}-(ii) we obtain

\begin{Lemma}\label{L:origine}
There exists $\rho\in (0,1)$ such that if  $u\in H^{1}(\R^{N})$ satisfies
$\norm{u}\leq \rho$
then $\f(u)\geq \tfrac 14\|u\|^2$ and
$\f'(u)v\geq \langle u,v\rangle-\tfrac 12\|u\|\|v\|$ for all $v\in H^{1}(\R^{N})$.
\end{Lemma}

\noindent By Lemma \ref{L:origine} in particular we get that
$\inf_{\|u\|=\rho}V(u)\geq\tfrac 14\rho^{2}>0$.
Moreover, by  Remark \ref{R:conseguenze}-(iii), we have that if $u\in H^1({\R^N})\setminus\{0\}$ and $s>1$  then
$V(su)=s^{2}\|u\|^{2}-\int_{R^{n}}F(su)\, dx\leq s^{2}\|u\|^{2}-s^{\mu}\int_{R^{n}}F(u)\, dx$ and hence $V(su)\to -\infty$ as $s\to+\infty$ for all $u\in H^1({\R^N})\setminus\{0\}$.
Hence, defining
\[\Gamma =\{\,\g\in C([0,1], H^1({\R^N}))\,:\ \g(0)=0\,,\ \gamma(1)\not=0\hbox{ and }\f(\g(1))\leq 0\,\}\] 
and setting
$$
c=\inf_{\g\in\Gamma}\max_{s\in[0,1]}\f(\g(s)),
$$
we get $c\geq \tfrac 14\rho^{2}$ and by the Mountain Pass Theorem (see e.g. \cite{[R]}) there exists a Palais Smale sequence for $\f$ at level $c$.\smallskip

\noindent 
Moreover, by $(f3)$, the following inequality holds true
\begin{equation}\label{eq:superV}
\mu V(u)-V'(u)u=(\tfrac\mu 2-1)\| u\|^{2}+\int_{\R^{N}}f(u)u-\mu F(u)\geq\tfrac{\mu-2}2\|u\|^{2},\end{equation}
from which in particular we derive that the Palais Smale sequences of $\f$ are bounded in $H^1({\R^N})$.  By (\ref{eq:superV}), we obtain also that if $\f'(u)=0$ and $u\not=0$ then $\f(u)\geq\tfrac{\mu -2}{2\mu}\| u\|^{2}$,
showing that $\f$ has not critical points (or Palais Smale sequences) at negative levels.
\smallskip

\noindent The existence of a mountain pass critical point of $V$ can then be deduced by using concentration compactness argument. We have
\begin{Proposition}\label{P:u0positive} There exists $w_{0}\in H^{1}(\R^{N})$ such that $\f(w_{0})=c$ and $\f'(w_{0})=0$. Moreover $w_{0}\in C^{2}(\R^{N})$ is a solution of (E) on $\R^{N}$,  $w_{0}>0$, $w_{0}(x)\to 0$ as $|x|\to+\infty$ and, up to translations, $w_{0}$ is radially symmetric about the origin with $\partial_{r}w_{0}<0$ for $r=|x|>0$.
\end{Proposition}
We refer for a proof to \cite{[BL]}, for $N\geq 3$ and \cite{[BGK]} for $N=2$, where a more general existence results regarding least energy solutions for scalar field equations is given. Their Mountain pass characterization is proved in \cite{[JJT]}. The case $N=1$ is easier and we omit the proof.\medskip\par\noindent
Fixed $u\in H^{1}(\R^{N})$, the assumption $(f4)$ allows us to describe the behaviour of $V$ along the rays $\{tu\,|\, t\geq 0\}$  
in $H^1({\R^N})$. 
\begin{Lemma}\label{L:raggi} For every $u\in H^1({\R^N})\setminus\{0\}$ there exists $t_{u}>0$ such that
\begin{equation}\label{eq:tu}\tfrac d{dt} \f(t u)>0\hbox{ for }t\in (0,t_{u})\hbox{ and }\tfrac d{dt} \f(tu)<0\hbox{ for }t\in (t_{u},+\infty).\end{equation}
Moreover $\f(t_{u}u)\geq c$ and for any $b\in (0,c)$ there exist unique  $\aub\in (0,t_{u})$ and $\oub\in (t_{u},
+\infty)$ such that $\f(\aub u)=\f(\oub u)=b$. Finally the function $t\mapsto \f'(tu)tu$ is decreasing in $(t_{u},+\infty)$.
\end{Lemma}
\Proof We have 
\begin{equation}\label{eq:derivataV(tu)}\tfrac d{dt} \f(t u)=\f'(tu)u=t(\|u\|^{2}-\tfrac 1t\int_{\R^{N}}f(tu)u\, dx).\end{equation}
By $(f4)$ 
the function $t\mapsto\tfrac 1t\int_{\R^{N}}f(tu)u\, dx$ is stricly increasing in $(0,+\infty)$ for any 
$u\not=0$ and so, by (\ref{eq:derivataV(tu)}), the fuction $\tfrac d{dt} \f(t u)$ can change sign at most in one point $t_{u}>0$. Then (\ref{eq:tu}) follows since $\f(0)=0$, $\f(su)\geq\tfrac 14 s^{2}\|u\|^{2}$ for $s\in (0,\rho/\|u\|)$ and $\f(su)\to-\infty$ as $s\to +\infty$.
By (\ref{eq:tu}) we deduce $\f(t_{u}u)=\max_{s\geq 0}\f(s u)$, and, 
by the definition of the mountain pass level, we have $\f(t_{u}u)\geq c$. Given $b\in [0,c)$, since $\f(0)=0$, $\f(tu)<0$ for $t$ large and $\f(t_{u}u)\geq c$, by continuity 
there exist (unique by (\ref{eq:tu})) $0\leq\aub<t_{u}<\oub$ such that $\f(\aub u)=\f(\oub u)=b$.
We finally note that by $(f4)$ we have $\tfrac{d^{2}}{dt^{2}}\f(tu)=\|u\|^{2}-\int_{\R^{N}}f'(tu)u^{2}\, dx\leq \|u\|^{2}-\tfrac 1t\int_{\R^{N}}f(tu)u\, dx<0$ for any $t>t_{u}$. 
We conclude that $\tfrac{d}{dt}\f'(tu)tu=\tfrac{d}{dt}(t\tfrac{d}{dt}\f(tu))=t\tfrac{d^{2}}{dt^{2}}\f(tu)+\tfrac{d}{dt}\f(tu)<0$ for any $t>t_{u}$.
\QED

\begin{Remark}\label{R:tucritici}{\rm Note that if $\f'(u)u=0$ and $u\not=0$ we have $\tfrac d{dt} \f(t u)\big|_{t=1}=\f'(u)u=0$ and so $t_{u}=1$. Then, by Lemma \ref{L:raggi}, $\f(u)=\f(t_{u}u)\geq c$ whenever $u\not=0$ and $\f'(u)u=0$. 
}\end{Remark}

\begin{Remark}\label{R:dimension} {\rm We note that, since by $(f2)$ we have $p<2^{*}_{N+1}-1$, all the results stated and proved in the present sections holds unchanged  for all $m\in\{1,\ldots,N+1\}$ considering the functionals
\[\f_{m}(u)=\tfrac12\|u\|_{H^{1}(\R^{m})}^2-\int_{{\R^m}} F(u(x))\,dx,\quad u\in H^{1}(\R^{m}).
\]
In particular, denoting $c_{m}$ the mountain pass level of $\f_{m}$ in $H^{1}(\R^{m})$, Proposition \ref{P:u0positive} establishes that $\f_{m}$ has a positive, radially symmetric, critical point $w\in H^{1}(\R^{m})$ at the level $c_{m}$.}\end{Remark}

\subsection{ Further properties of $\f$ on the space of radial functions. The sublevels $\VV^{b}_{-}$ and $\VV^{b}_{+}$}
From now on we reduce ourself to work on the subspace of $H^{1}$ constituted by radial functions: $H^1_{r}({\R^N})=\{u\in H^{1}(\R^{N})\,/\, u(x)=u(|x|)\}$.  We recall that by the Strauss Lemma (see \cite{[S]}, \cite{[L]}) $H^1_{r}({\R^N})$  is compactly embedded in $L^{q}(\R^{N})$ for all $q\in (2, 2^{*}_{N})$. Thanks to the Strauss Lemma the functional $\f$ is weakly lower semicontinuous on $H^1_{r}({\R^N})$.
\begin{Lemma}\label{L:semiV}
 Let $u_{n}\to u$ and $v_{n}\to v$ weakly in $H^1_{r}({\R^N})$. Then 
 \[
 \lim_{n\to+\infty}\int_{\R^{N}} F(u_{n})dx=\int_{\R^{N}} F(u)dx\quad\hbox{and}\quad\lim_{n\to+\infty}\int_{\R^{N}} f(u_{n})v_{n}\, dx=\int_{\R^{N}} f(u)vdx.\]
Hence $\f(u)\leq{\displaystyle\liminf_{n\to+\infty}}\,\f(u_{n})$, $\f'(u)u\leq{\displaystyle\liminf_{n\to+\infty}}\,\f'(u_{n})u_{n}$ and, for every $h\in H^1_{r}({\R^N})$, $\f'(u)h={\displaystyle\lim_{n\to+\infty}}\,\f'(u_{n})h$.
\end{Lemma}
\Proof Since, by $(f2)$, $p+1<2^{*}_{N}$, we have $H^1_{r}({\R^N})$ is compactly embedded in $L^{p+1}(\R^{N})$ and so $u_{n}\to u$ and $v_{n}\to v$ strongly in $L^{p+1}(\R^{N})$. Then, since by (\ref{eq:originf}) we have that for all $\e>0$
\begin{align*} |F(u_{n})-F(u)|&= |\int_0^1
f(u+s(u_{n}-u))(u_{n}-u)\, ds|\\
&\leq \epsilon| u_{n}-u| (| u| +|u_{n}-u| )+2^{p-1}A_{\epsilon}
| u_{n}-u|( | u|^p+| u_{n}-u|^p),\end{align*}
we deduce that  as $n\to+\infty$
\begin{align*}
&\int_{\R^{N}} |F(u_{n})-F(u)|\, dx
\leq\epsilon\|u_{n}-u\|_{2}(\|u\|_{2}+\|u_{n}-u\|_{2})\\
&\phantom{\int_{\R^{N}} }+
2^{p-1}A_{\epsilon}\|u_{n}-u\|_{p+1}(\|u\|_{p+1}^{p}+\|u_{n}-u\|_{p+1}^{p})
=\epsilon\|u_{n}-u\|_{2}(\|u\|_{2}+\|u_{n}-u\|_{2})+o(1).\end{align*}
Since $\epsilon$ is arbitrary and $\|u_{n}-u\|$ bounded, we deduce that $\int_{\R^{N}} F(u_{n})\, dx\to\int_{\R^{N}} F(u)dx$. \\
We show now that $\int_{\R^{N}} f(u_{n})v_{n}\, dx\to\int_{\R^{N}} f(u)vdx$.  
First we write 
\begin{equation}\label{eq:decomp}\int_{\R^{N}} f(u_{n})v_{n}\, dx=\int_{\R^{N}} f(u_{n})(v_{n}-v)dx+\int_{\R^{N}} f(u_{n})v dx.\end{equation} 
We note that by (\ref{eq:originf}), since $\|v_{n}-v\|_{p+1}\to 0,$ as $n\to+\infty$ we have 
\begin{align*}|\int_{\R^{N}} f(u_{n})(v_{n}-v)& dx|\leq\int_{\R^{N}}\epsilon |u_{n}||v_{n}-v|+A_{\epsilon}|u_{n}|^{p}|v_{n}-v|\, dx\\
&\leq\epsilon\|u_{n}\|_{2}\|v_{n}-v\|_{2}+A_{\epsilon}\|u_{n}\|_{p+1}^{p}\|v_{n}-v\|_{p+1}=\epsilon\|u_{n}\|_{2}\|v_{n}-v\|_{2}+o(1)\end{align*}
and so, since $\epsilon$ is arbitrary, we deduce $\int_{\R^{N}} f(u_{n})(v_{n}-v) dx\to 0$.
Then, by (\ref{eq:decomp}), our claim follows if we show that $\int_{\R^{N}} f(u_{n})v dx\to \int_{\R^{N}} f(u)v dx$.
For that
we fix $\epsilon>0$ and choose $R>0$ such that $\int_{|x|>R}|v|^{2}\, dx+\int_{|x|>R}|v|^{p+1}\, dx<\epsilon$. By (\ref{eq:originf}), with $\varepsilon=1$, we have
\begin{align*}|\int_{|x|>R}f(u_{n})v-&f(u)v\, dx|\leq \int_{|x|>R}(|u_{n}|+|u|)|v|+A_{1}( | u|^p+| u_{n}|^p)|v|\, dx\\
&\leq \epsilon^{1/2}(\|u_{n}\|_{2}+\|u\|_{2})+A_{1}\epsilon^{1/(p+1)}(\|u_{n}\|^{p}_{p+1}+\|u\|^{p}_{p+1})\end{align*}
 and, since $\epsilon$ is arbitrary,  we deduce $\int_{|x|>R}f(u_{n})v\, dx\to \int_{|x|>R}f(u)v\, dx$.
On the other hand we have $u_{n}\to u$ in $L^{2}(B_{R}(0))$ and in $L^{p+1}(B_{R}(0))$. Then for any subsequence $(u_{n_{k}})$ there exists
 a subsubsequence $(u_{n_{k_{j}}})$  and  a function $\psi\in L^2(B_{R}(0))\cap
L^{p+1}(B_{R}(0))$ such that
$u_{n_{k_j}}(x)\to u(x)$ a.e. in $B_{R}(0)$ and $|
u_{n_{k_j}}(x)|\leq \psi(x)$ a.e. in $B_{R}(0)$, for any $j\in\N$.
Using again (\ref{eq:originf}) we obtain $|f(u_{n_{k_{j}}})-f(u)||v|\leq (|\psi|+|u|)|v|+A_{1}( | \psi|^p+|u|^p)|v|$ on $B_{R}(0)$ for any $j\in\N$. Then, by
the dominated convergence theorem, we get $\int_{|x|<R}f(u_{n_{k_{j}}})v\, dx\to \int_{|x|<R}f(u)v\, dx$. Since the subsequence $(u_{n_{k}})$ is arbitrary we conclude that $\int_{|x|<R}f(u_{n})v\, dx\to \int_{|x|<R}f(u)v\, dx$
and the Lemma follows.
 \QED\bigskip

For our study it is important to understand the structure of the sublevel sets $\VV^{b}=\{u\in H^1_{r}({\R^N})\,/\, \f(u)\leq b\}$. By definition of the Mountain pass level the set $\VV^{b}$ is not path connected for any $b\in [0,c)$. Given $b\in [0,c)$, recalling Lemma \ref{L:raggi}, 
 we denote 
\[\VV_{-}^{b}=\{ tu\,|\, u\in H^1_{r}({\R^N})\setminus\{0\},\, t\in[0, \aub] \}\hbox{ and } \VV_{+}^{b}=\{ tu\,|\, u\in H^1_{r}({\R^N})\setminus\{0\},\, t\in[\oub, +\infty)\}.\]
Clearly 
\[\VV^{b}=\VV_{-}^{b}\cup\VV_{+}^{b}.\]
\begin{Remark}\label{R:connected} {\rm The set $\VV_{-}^{b}$ is clearly path connected (starshaped indeed, with respect to the origin). The same holds true also for $\VV_{+}^{b}$. Indeed, given $u_{1}, u_{2}\in\VV_{+}^{b}$ such that $b\geq b_{1}=\f(u_{1})\geq b_{2}=\f(u_{2})$
we can connect them considering the path $\gamma(s)=\omega_{b_{1},(1-s)u_{1}+su_{2}}((1-s)u_{1}+su_{2})$ for $s\in [0,1]$ and
$\gamma(s)=s\omega_{b_{1}.u_{2}}$ for $s\in [1,1/\omega_{b_{1},u_{2}}]$. The function $\gamma$ is continuous since
the mapping $u\in H^1_{r}({\R^N})\to\oub\in\R$ is continuous for any $b< c$.}\end{Remark}
\begin{Remark}\label{R:caratterizzazioneVb}{\rm By definition of mountain pass level and Remark \ref{R:connected}, if $\gamma\in C([0,1],H^1_{r}({\R^N}))$ is such that $\gamma(0)\in \VV_{-}^{b}$ and $\gamma(1)\in\VV_{+}^{b}$ then $\max_{s\in[0,1]}\f(\gamma(s))\geq c$. Secondly note that by Lemma \ref{L:raggi}
\[\VV_{-}^{b}=\{u\in H^1_{r}({\R^N})\,/\, \aub\geq 1\}\cup\{ 0\}\hbox{ and }\VV_{+}^{b}=\{u\in H^1_{r}({\R^N})\,/\, \oub\leq 1\}\hbox{ for all }b\in [0,c).\]
Moreover  if $b\in (0,c)$ then
\begin{equation}\label{eq:V-}u\in \VV_{-}^{b}\setminus\{0\}\hbox{ if and only if }\f(u)\leq b\hbox{ and }\f'(u)u> 0.\end{equation}
Indeed,  if $u\in \VV_{-}^{b}\setminus\{0\}$ then $1\leq\aub<t_{u}$ and so, by Lemma \ref{L:raggi},  $\f'(u)u>0$. Viceversa
if
$\f(u)\leq b$ and $\f'(u)u> 0$ then $u\not=0$ and $1\leq\aub$, from which $V(u)\le b$.
Analogously if $b\in [0,c)$ then
\begin{equation}\label{eq:Vp}u\in \VV_{+}^{b}\hbox{ if and only if }\f(u)\leq b\hbox{ and }\f'(u)u< 0.\end{equation}}\end{Remark}
\begin{Lemma}\label{L:chisuradebole}
If $b\in [0,c)$ then $\VV_{-}^{b}$ and $\VV_{+}^{b}$ are weakly closed in $H^1_{r}({\R^N})$.\end{Lemma}
\Proof
Let $(u_{n})\subset\VV^{b}_{+}$ be such that $u_{n}\to u_{0}$ weakly in $H^1_{r}({\R^N})$. By Remark \ref{R:caratterizzazioneVb} we have $\f(u_{n})\leq b$ and $\f'(u_{n})u_{n}< 0$. Since $\f'(u_{n})u_{n}< 0$, by Lemma \ref{L:origine} we deduce $\|u_{n}\|\geq\rho$ for any $n\in\N$. Moreover since
$\f(u_{n})\leq b$, by Lemma  \ref{L:semiV} we obtain $\f(u_{0})\leq b$.
By Lemma \ref{L:semiV} we know also that $\int_{\R^{n}}f(u_{n})u_{n}\, dx\to\int_{\R^{n}}f(u_{0})u_{0}\, dx$ and, since $\f'(u_{n})u_{n}< 0$,  $
\f'(u_{0})u_{0}\leq 0$. By (\ref{eq:Vp}), to prove that $u_{0}\in\VV^{b}_{+}$ we have to show that $\f'(u_{0})u_{0}<0$. For that, assume by contradiction that $\f'(u_{0})u_{0}=0$ and note that, being $\f(u_{0})\leq b<c$, by Remark \ref{R:tucritici} we have $u_{0}=0$. Then $\int_{\R^{n}}f(u_{n})u_{n}\, dx\to 0$ and so $0>\f'(u_{n})u_{n}>\rho^{2}+o(1)$ as $n\to+\infty$, a contradiction which shows that $\VV^{b}_{+}$ is weakly closed.\\
Let now $(u_{n})\subset\VV^{b}_{-}$ be such that $u_{n}\to u_{0}$ weakly in $H^1_{r}({\R^N})$. Again using Remark \ref{R:caratterizzazioneVb} we have $\f(u_{n})\leq b$ and $\f'(u_{n})u_{n}\geq 0$. Hence, by Lemma \ref{L:semiV}, we deduce that $\f(u_{0})\leq b$. To show that $u_{0}\in\VV_{-}^{b}$ it suffices to show that $\f'(u_{0})u_{0}\geq 0$.
Assume by contradiction that $\f'(u_{0})u_{0}< 0$. Then, by (\ref{eq:Vp}),
 we have $u_{0}\in\VV_{+}^{b}$. Consider the path $\gamma_{n}(s)=u_{0}+s(u_{n}-u_{0})$, $s\in [0,1]$.
Since $\gamma_{n}(0)=u_{0}\in\VV_{+}^{b}$ and  $\gamma_{n}(1)=u_{n}\in\VV_{-}^{b}$, by Remark \ref{R:caratterizzazioneVb},
for any $n\in\N$ we find $s_{n}\in (0,1)$ such that $\f(\gamma_{n}(s_{n}))\geq c$. 
We note also that
$\|\gamma_{n}(s)\|_{2}\leq\|u_{0}\|_{2}+\|u_{n}-u_{0}\|_{2}\leq C_{1}<+\infty$  and
$\|\gamma_{n}(s)\|_{p+1}\leq\|u_{0}\|_{p+1}+\|u_{n}-u_{0}\|_{p+1}\leq C_{2}<+\infty$ for any $n\in\N$ and $s\in [0,1]$. Then, choosing $\epsilon=\tfrac{c-b}{2C_{1}^{2}}$, by (\ref{eq:originf}) we get
\begin{align*}|\int_{\R^{N}} f(\gamma_{n}(s))(u_{n}-u_{0}) dx|&\leq\epsilon
\|\gamma_{n}(s)\|_{2}\|u_{n}-u_{0}\|_{2}+A_{\epsilon}\|\gamma_{n}(s)\|_{p+1}^{p}\|u_{n}-u_{0}\|_{p+1}\\
&=\tfrac{c-b}{2}+A_{\epsilon}C_{2}^{p}\|u_{n}-u_{0}\|_{p+1}\quad\hbox{ for any }s\in [0,1].\end{align*}
Hence we derive that for any $s\in[0,1]$ and $n\in\N$ there results
\begin{align*}
\tfrac{d}{ds}\f(\gamma_{n}(s))&=\f'(\gamma_{n}(s))(u_{n}-u_{0})\\&\geq
s\|u_{n}-u_{0}\|^{2}+\langle u_{0}, u_{n}-u_{0}\rangle-\tfrac{c-b}{2}-A_{\epsilon}C_{2}^{p}\|u_{n}-u_{0}\|_{p+1}.\end{align*}
Integrating on $[s_{n},1]$ we get
\begin{align*}b-c&\geq\f(u_{n})-\f(\gamma_{n}(s_{n}))\geq \tfrac{b-c}{2}+(1-s_{n})(\langle u_{0}, u_{n}-u_{0}\rangle-
A_{\epsilon}C_{2}^{p}\|u_{n}-u_{0}\|_{p+1}).\end{align*}
Since $\langle u_{0}, u_{n}-u_{0}\rangle-A_{\epsilon}C_{2}^{p}\|u_{n}-u_{0}\|_{p+1} \to 0$ we obtain the contradiction $0>b-c\geq \tfrac{b-c}{2}.$
\QED
\begin{Remark}\label{R:V-bound}{\rm Note that, by (\ref{eq:superV}),  if $b\in[0,c)$ and $u\in \VV^{b}_{-}$, since $\f'(u)u\geq 0$, then
$$
\|u\|^{2}\leq \tfrac{2\mu}{\mu-2}V(u)\le\tfrac{2\mu}{\mu-2}b.
$$
In particular we obtain that $\VV^{b}_{-}$ is bounded in $H^1_{r}({\R^N})$. Then, by Lemma \ref{L:chisuradebole}, $\VV^{b}_{-}$ is weakly compact in $H^1_{r}({\R^N})$ and if $(u_{n})\subset\VV^{b}_{-}$ is such that
$u_{n}\to u_{0}$ with respect to the $L^{2}(\R^{N})$ metric then $u_{0}\in \VV^{b}_{-}$.}\end{Remark}

\begin{Lemma}\label{L:gradienteV+} If $b\in [0,c)$ we have 
$\nu^{+}(b):=\inf_{u\in\VV^{b}_{+}}\tfrac{-\f'(u)u}{\max\{1,\|u\|^{2}_{2}\}}>0$. 
\end{Lemma}
\Proof First note that, by (\ref{eq:superV}), if $u\in\VV^{b}_{+}$ is such that $\|u\|_{2}^{2}\geq\tfrac{4b\mu}{\mu-2}$ or $\f(u)\leq 0$ then
$\tfrac{-\f'(u)u}{\|u\|^{2}_{2}}\geq\tfrac{\mu-2}2\tfrac{\|u\|^{2}}{\|u\|^{2}_{2}}-\mu \tfrac{\f(u)}{\|u\|_{2}^{2}}\geq
\tfrac{\mu-2}4$.
Assume now by contradiction that there exists $(u_{n})\subset\VV^{b}_{+}$ such that $0<\f(u_{n})\leq b$, $ \|u_{n}\|_{2}^{2}\leq\tfrac{4b\mu}{\mu-2}$ and $\tfrac{-\f'(u_{n})u_{n}}{\max\{1,\|u_{n}\|^{2}_{2}\}}\to 0$.  Then $\f'(u_{n})u_{n}\to 0$. Since $u_{n}\in\VV^{b}_{+}$, by Remark \ref{R:caratterizzazioneVb} we have $t_{u_{n}}<1$. By (\ref{eq:superV}) we have
$\|u_{n}\|^{2}\leq \tfrac{2\mu}{\mu -2}b+o(1)$ and since, by Remark \ref{R:tucritici}, 
$\|t_{u_{n}}u_{n}\|\geq\rho$, we deduce that $t_{u_{n}}\geq \tfrac{\mu -2}{4\mu b}\rho>0$ whenever $n$ is large. By Lemma \ref{L:raggi} we have
$|\f'(su_{n})su_{n}|\leq |\f'(u_{n})u_{n}|$ for any $s\in (t_{u_{n}},1)$,  and we conclude
$c-b\leq\int_{1}^{t_{u_{n}}}\tfrac d{ds}\f(su_{n})\, ds=\int_{1}^{t_{u_{n}}}\tfrac 1s\f'(su_{n})su_{n}\, ds\leq -\log(t_{u_{n}})|\f'(u_{n})u_{n}|\to 0\hbox{ as }n\to+\infty$,
a contradiction which proves the Lemma.\QED

\begin{Lemma}\label{L:gradienteV-} If $b\in (0,c)$ then
$\nu^{-}(b):=\inf_{u\in\VV^{(b+c)/2}_{-}\setminus \VV^{ b}_{-}}\f'(u)u
>0$. 
\end{Lemma}
\Proof 
By contradiction, let $(u_{n})\subset\VV^{(b+c)/2}_{-}\setminus \VV^{ b}_{-}$ be such that  $\f'(u_{n})u_{n}\to 0$. Then,  
by Remark \ref{R:V-bound}, there exists $u_{0}\in \VV^{(b+c)/2}_{-}$ such that, up to a subsequence, $u_{n}\to u_{0}$ weakly in $H^1_{r}({\R^N})$. By Lemma \ref{L:semiV},  $\f'(u_{0})u_{0}\le\liminf \f'(u_{n})u_{n}=0 $. Since $u_{0}\in \VV^{(b+c)/2}_{-}$ that implies $u_{0}=0$ and
then, again by Lemma \ref{L:semiV}, $\int_{\R^{n}}f(u_{n})u_{n}\, dx\to 0$. Hence $\f'(u_{n})u_{n}=\|u_{n}\|^{2}+o(1)\to 0$ and so
$u_{n}\to 0$ in $H^{1}(\R^{n})$ that gives the contradiction
$0<b\leq\f(u_{n})\to 0$.\QED
\medskip

Finally, we display some properties depending on the assumption $p<1+\frac 4N$.\bigskip

\noindent First, as a particular case of the Gagliardo Nirenberg interpolation inequality (see \cite{[NiGa]}), we have that there exists a constant $\kappa=\kappa(N,p)>0$ such that for any $u\in H_{r}^{1}(\R^{N})$, there results
\begin{equation}\label{eq:interpolazione}\|u\|_{p+1}\leq\kappa\|u\|_{2}^{\theta}\|\nabla u\|_{2}^{1-\theta},\quad\hbox{where }1-\theta=\tfrac N2\tfrac {p-1}{p+1}.\end{equation}
Moreover, note that, by (\ref{eq:originF}), we have 
$F(t)\leq \tfrac 14 |t|^{2}+\tfrac{A_{1/2}}{p+1}|t|^{p+1}$ for every $t\in\R$.
Therefore, if $u\in H^1_{r}({\R^N})\setminus\{ 0\}$, by (\ref{eq:interpolazione}) there results
\begin{equation}\label{eq:Vmaggiore}\f(u)\geq \tfrac 12\|\nabla u\|_{2}^{2} \left(1-\tfrac{2\kappa_{GN}A_{1/2}}{p+1}\tfrac{\| u\|_{2}^{(p+1)\theta}}{\|\nabla u\|_{2}^{2-(p+1)(1-\theta)}}\right)+
\tfrac 14\|u\|_{2}^{2}.\end{equation}
where, since $p<1+\frac 4N$, by $(f2)$, we have
\begin{equation}\label{eq:1menthetap+1}
(p+1)(1-\theta)=\tfrac N2 (p-1)<2.
\end{equation}
By (\ref{eq:Vmaggiore}) and (\ref{eq:1menthetap+1}) it follows directly
\begin{Lemma}\label{L:coercivita}
If $(u_{n})\subset H^1_{r}({\R^N})$,  
$\displaystyle\sup_{n\in\N}\| u_{n}\|_{2}<+\infty$ and  $\| \nabla u_{n}\|_{2}\to+\infty$ 
then $\f(u_{n})\to+\infty$.
\end{Lemma}
In particular $\VV^{b}_{+}$ enjoys the following property.
\begin{Lemma}\label{L:bound} 
If $b\in [0,c)$, for any $M_{1}>0$ there exists $M_{2}>0$ such that
if $u\subset \VV^{b}_{+}$ and $\|u\|_{2}\leq M_{1}$ then $\|\nabla u\|_{2}\leq M_{2}$.\end{Lemma}
\begin{Remark}\label{R:V-bound+}{\rm Note that by Lemma \ref{L:bound} and Lemma
\ref{L:chisuradebole} we derive that  if $(u_{n})\subset\VV^{b}_{+}$ is such that
$u_{n}\to u_{0}$ with respect to the $L^{2}(\R^{N})$ metric then $u_{0}\in \VV^{b}_{+}$.}\end{Remark}
Another consequence is the following one
\begin{Lemma}\label{L:distanzaL2} 
For any $b_{1},\, b_{2}\in[0,c)$ there result $\delta(b_{1},b_{2}):=\mathrm{dist}(\VV_{-}^{b_{1}},\VV_{+}^{b_{2}})>0$.
\end{Lemma}
\Proof Clearly $\delta(b_{1},b_{2})<+\infty$. Let $(u_{n,1})\subset\VV_{-}^{b_{1}}$ and $(u_{n,2})\subset\VV_{+}^{b_{2}}$ be such that
$\|u_{n,1}-u_{n,2}\|_{2}\to \delta(b_{1},b_{2})$. By Remark \ref{R:V-bound} we know that
$\|u_{n,1}\|\leq \tfrac{2\mu}{\mu -2}b_{1}$ and hence we obtain
$\|u_{n,2}\|_{2}\leq  \tfrac{2\mu}{\mu -2}b_{1}+\delta(b_{1},b_{2})+o(1)$. Then $(u_{n,2})$ is bounded in $L^{2}(\R)$.
By Lemma \ref{L:coercivita}, since $\f(u_{n,2})\leq b_{2}$, we obtain that
$\sup_{n\in\N}\|\nabla u_{n,2}\|_{2}<+\infty$ and so that $(u_{n,2})$ is bounded also in $H^1_{r}({\R^N})$. Then there exists two subsequences $(u_{n_{j},1})\subset (u_{n,1})$, $(u_{n_{j},2})\subset (u_{n,2})$ which weakly converge respectively to $u_{1}\in H^1_{r}({\R^N})$ and $u_{2}\in H^1_{r}({\R^N})$. By Lemma
 \ref{L:chisuradebole} we have $u_{1}\in \VV^{b_{1}}_{-}$ and $u_{2}\in \VV^{b_{2}}_{+}$ and by the weak semicontinuity of the $L^{2}$ norm we deduce
 $ \delta(b_{1},b_{2})\leq\|u_{1}-u_{2}\|_{2}\leq\lim_{j\to+\infty}\|u_{n_{j},1}-u_{n_{j},2}\|_{2}=\delta(b_{1},b_{2})$.
Since $u_{1}\not=u_{2}$ we have $\delta(b_{1},b_{2})=\|u_{1}-u_{2}\|_{2}>0$ and the Lemma follows.\QED

As a further consequence of the assumption $p<1+4/N$, we give a result  concerning the behaviour of $\f$ along sequences in $H^1_{r}({\R^N})$ which converge to a point $u_{0}\in H^1_{r}({\R^N})$ with respect to the $L^{2}(\R^{N})$ metric.
\begin{Lemma}\label{L:taglio-conv-L2}
Let $u_{n},u_{0} \in H^1_{r}({\R^N})$ be such that 
$\|u_{n}- u_{0}\|_{2}\to 0$  as $n\to+\infty$
 and $\liminf_{n\to\infty}\|\nabla (u_{n}- u_{0})\|_{2}>0$. 
 Then 
there exists $\bar n\in\N$ such that  
$$
    \f(u_{n})-\f(u_{0}+s(u_{n}-u_{0}))\geq \tfrac14 (1-s)\|\nabla (u_{n}-u_{0})\|_{2}^{2},
    \qquad \fa s\in [0,1],\, n\geq\bar n.
$$
\end{Lemma}
\Proof
\noindent Setting $w_{n}=u_{n}-u_{0}$,  by (\ref{eq:originf}), since $w_{n}\to 0$ in $L^{2}(\R^{N})$ we recover that there exists $C>0$ such that, for any $s\in [0,1]$,
\begin{align}\nonumber
   &|\int_{\R^{N}}F(u_{0}+w_{n})-F(u_{0}+sw_{n}))\,
   dx|=|\int_{\R^{N}}\int_{s}^{1}f(u_{0}+\sigma w_{n})w_{n} d\sigma\,dx|\\ \nonumber
   &\leq \int_{s}^{1}\|u_{0}\|_{2}\|w_{n}\|_{2}+\sigma\|w_{n}\|_{2}^{2}+
   A_{1}2^{p-1}(\|u_{0}\|^{p}_{p+1}\|w_{n}\|_{p+1}+\sigma^{p}\|w_{n}\|^{p+1}_{p+1})\,d\sigma\\\label{eq:Wto0}
   &\leq C(1-s)(o(1)+\|w_{n}\|_{p+1}+\|w_{n}\|^{p+1}_{p+1})\qquad\hbox{ as }n\to+\infty.
\end{align}
We now note that, since  $\liminf_{n\to+\infty}\|\nabla w_{n}\|_{2}^{2}>0$,  we have
\begin{equation}\label{eq:stimegradscalar}
   \lim_{n\to+\infty}\tfrac{\langle\nabla u_{0},\nabla w_{n}\rangle_{2}}{\|\nabla w_{n}\|_{2}^{2}} =0.\end{equation} 
Indeed, (\ref{eq:stimegradscalar}) is true along subsequences $(w_{n_{j}})$ such that $\|\nabla w_{n_{j}}\|_{2}\to+\infty$.
If $(w_{n_{j}})\subset\{w_{n}\}$ is bounded in $H^1_{r}({\R^N})$ then, necessarily, $w_{n_{j}}\to 0$ weakly in $H^1_{r}({\R^N})$ and again  (\ref{eq:stimegradscalar})  follows.\smallskip

\noindent Secondly we note that
\begin{equation}\label{eq:rapporto}
\lim_{n\to+\infty}\tfrac{\|w_{n}\|_{p+1}+\|w_{n}\|^{p+1}_{p+1}}{\|\nabla w_{n}\|_{2}^{2}}= 0.\end{equation}
Indeed,  we have either $\|\nabla w_{n}\|_{2}$ is bounded or $\limsup_{n\to+\infty}\|\nabla w_{n}\|_{2}= +\infty$.
If $\|\nabla w_{n}\|_{2}$ is bounded then $(w_{n})$ weakly converges to $0$ in $H^1_{r}({\R^N})$ and so strongly in $L^{p+1}(\R^{N})$ giving (\ref{eq:rapporto}). If $\|\nabla w_{n}\|_{2}\to +\infty$ along a subsequence, then, since $\|w_{n}\|_{2}\to 0$, (\ref{eq:rapporto}) follows by (\ref{eq:interpolazione}) and (\ref{eq:1menthetap+1}).\smallskip

\noindent Finally, by (\ref{eq:Wto0}), we derive that for any 
$s\in [0,1]$
\begin{align*}
    &\f(u_{0}+w_{n})-\f(u_{0}+sw_{n})=
   \tfrac{\|\nabla w_{n}\|^{2}}{2}(1-s^{2})+
    (1-s)\langle\nabla u_{0},\nabla w_{n}\rangle_{2}+(1-s)o(1)
    \\
    &\phantom{\f(u_{0}+w_{n})-\f(u_{0}+sw_{n})=
   \tfrac{\|\nabla w_{n}\|^{2}}{2}(1-s^{2})}
   -\int_{\R^{N}}F(u_{0}+w_{n})-F(u_{0}+sw_{n})\,dx\\
    &\geq \|\nabla w_{n}\|_{2}^{2}(1-s)(\tfrac{1+s}{2}+\tfrac{\langle\nabla u_{0},\nabla w_{n}\rangle_{2}}{\|\nabla w_{n}\|_{2}^{2}}-C
    \tfrac{\|w_{n}\|_{p+1}+\|w_{n}\|^{p+1}_{p+1}+o(1)}{\|\nabla w_{n}\|_{2}^{2}})\geq \|\nabla w_{n}\|_{2}^{2}(1-s)(\frac12+o(1))
\end{align*}
and the Lemma follows by (\ref{eq:stimegradscalar}) and (\ref{eq:rapporto}).
\QED

\begin{Remark}\label{R:conv-H1}{\rm By Lemma \ref{L:taglio-conv-L2} we have in particular that
if $u_{n},u_{0} \in H^1_{r}({\R^N})$, $s_{n}\in[0,1]$ are such that $u_{n}\to u_{0}$ in $L^{2}(\R^{N})$ as $n\to+\infty$
and
$\f(u_{n})- \f(u_{0}+s_{n}(u_{n}-u_{0}))\to 0$ as $n\to+\infty$,
then $(1-s_{n})\| u_{n}- u_{0}\|^{2}\to 0$ as $n\to +\infty$.
In particular, if $\f(u_{n})\to \f(u_{0})$ as $n\to+\infty$,
then $u_{n}\to u_{0}$ in $H^{1}(\R^{N}$ as $n\to +\infty$.}
\end{Remark}

\section{Solutions on $\R^{N+1}$.}

In the sequel we denote $(x,y)\in\R^{N+1}$ where
$x=(x_{1},\ldots,x_{n})\in\R^{N}$ and $y\in\R$, the gradient with respect to the $x\in\R^{N}$ will be denoted by $\nabla_{x}$.
For $(y_{1},y_{2})\subset\R$ we set
$S_{(y_{1},y_{2})}:=\R^{N}\times (y_{1},y_{2})$ and, more simply, $S_{L}:=S_{[-L,L]}$ for $L>0$. We denote by $\X$ the set of monotone decreasing radially symmetric functions in $H^{1}(\R^{N})$:
$$
\X=\{u\in H^{1}_{r}(\R^{N})\,|\, u(x_{1})\geq u(x_{2})\hbox{ for any }x_{1},\, x_{2}\in\R^{N}\hbox{ such that }|x_{1}|\leq |x_{2}|\}.
$$
Note that $\X$ is a positive cone in $H^{1}_{r}(\R^{N})$ (and so convex) and it is sequentially closed in $H^{1}(\R^{N})$ with respect to the weak topology. In the following, with abuse of notation, given $b\in [0,c)$ we will indicate $\VV^{b}_{\pm}\equiv \VV^{b}_{\pm}\cap \X$.
\medskip

We consider the set
\[\MM=\{ v\in \cap_{L>0}H^{1}(S_{L})\,/\, v(\cdot,y)\in \X\hbox{ for a.e. }y\in\R\}.
\]
Note that, by Fubini Theorem, we have
that if $v\in \MM$
then 
$v(x,\cdot)\in H^{1}_{loc}(\R)$ for a.e. $x\in\R^{N}$. Therefore,
if $(y_{1},y_{2})\subset\R$ then
$v(x,y_{2})-v(x,y_{1})=\int_{y_{1}}^{y_{2}}\partial_{y}v(x,y)\, dy$
holds for a.e. $x\in\R^{N}$ and so
\[
\int_{\R^{N}}|v(x,y_{2})-v(x, y_{1})|^{2}\, dx=\int_{\R^{N}} |\int_{y_{1}}^{y_{2}}\partial_{y}u(x,y)\, dy\,|^{2} dx
	\leq 
	|y_{2}-y_{1}|\int_{\R^{N}}
	\int_{y_{1}}^{y_{2}}|\partial_{y}v(x,y)|^{2}\,dydx
\]
According to that, if $v\in
\MM$, the function $y\in\R\mapsto u(\cdot,y)\in L^{2}(\R^{N})$, defines a continuous 
trajectory
verifying
\begin{equation}\label{eq:continuita1}
    \|v(\cdot, y_{2})-v(\cdot, y_{1})\|_{2}^{2} \leq
    \|\partial_{y} v\|^{2}_{L^{2}
    (S{(y_{1},y_{2})})}
    |y_{2}-y_{1}|,\quad\forall\, (y_{1},y_{2})\subset\R.
\end{equation}
In the sequel we will consider the functional $\f$ as extended on $L^{2}(\R^{N})$ in the following way
\[
\f(u)=\begin{cases}
    \f(u), &\hbox{if } u\in H^{1}(\R^{n}), \\
    +\infty,  &\hbox{if } u\in L^{2}(\R^{N})\setminus H^{1}(\R^{n}).
    \end{cases}
\]

\begin{Lemma}\label{L:semicVv}
If $v\in\MM$ then the function $y\in\R\to\f(v(\cdot,y))\in \R\cup\{+\infty\}$ is lower semicontinuous.\end{Lemma}
\Proof
Let $v\in\MM$ and $y_{n}\to y_{0}$ and let $(y_{n_{j}})\subset(y_{n})$ be such that 
$\liminf_{n\to+\infty}\f(v(\cdot,y_{n}))=\lim_{j\to+\infty}\f(v(\cdot,y_{n_{j}}))
$. 
By (\ref{eq:continuita1}) we have $v(\cdot,y_{n_{j}})\to v(\cdot,y_{0})$ in $L^{2}(\R^{N})$  as $j\to+\infty$. We consider the two following alternative case: 
\[ \hbox{ (a) }\sup_{j\in\N}\|v(\cdot,y_{n_{j}})\|<+\infty\quad\hbox{ or }\quad \hbox{(b) }\limsup_{j\to +\infty}\|v(\cdot,y_{n_{j}})\|=+\infty\]
In the case (a), since
$(v(\cdot,y_{n_{j}}))$ is bounded in $\X$ and $v(\cdot,y_{n_{j}})\to v(\cdot,y_{0})$ in $L^{2}(\R^{N})$, we deduce that $v(\cdot,y_{n_{j}})\to v(\cdot,y_{0})$ weakly in $\X$. Then by Lemma
\ref{L:semiV} we derive
$\lim_{j\to+\infty}\f(v(\cdot,y_{n_{j}}))\geq\f(v(\cdot,y_{0}))$. In the case (b) we have $\limsup_{j\to+\infty}\|\nabla v(\cdot,y_{n_{j}})\|_{2}=+\infty$ since $\|v(\cdot,y_{n_{j}})\|_{2}$ is bounded. Then, by Lemma \ref{L:coercivita}, we get
$\lim_{j\to+\infty}\f(v(\cdot,y_{n_{j}}))=\limsup_{j\to+\infty}\f(v(\cdot,y_{n_{j}}))=+\infty$,
showing that also in the case (b) there results $\lim_{j\to+\infty}\f(v(\cdot,y_{n_{j}}))\geq \f(v(\cdot,y_{0}))$.
\QED

\begin{Lemma}\label{L:energia} If $v\in \MM$ is a solution of
 (E) on $S_{(y_{1},y_{2})}$ then
the energy function $E_{v}(y)=\tfrac{1}{2}\|\partial_{y}v(\cdot,y)\|_{2}^{2}-\f(v(\cdot,y))$ is constant on
$(y_{1},y_{2})$.\end{Lemma}
\Proof
Since $v\in\MM$ we have that $v\in H^{1}(S_{L})$ for any $L>0$. Then, since
$v$ solve $(E)$ on $S_{(y_{1},y_{2})}$ by regularity we have $v\in H^{2}(S_{(\zeta_{1},\zeta_{2})})\cap
C^{2}(S_{(y_{1},y_{2})})$ for any $[\zeta_{1},\zeta_{2}]\subset (y_{1},y_{2})$. Hence $v(\cdot,y)\in H^{2}(\R^{N})\cap C^{2}(\R^{N})$ for all $y\in (y_{1},y_{2})$  
and so
$\int_{|x|=R}|v(x,y)|+|\nabla_{x}v(x,y)|\, d\sigma\to 0$ as ${R\to+\infty}$ for all $y\in (y_{1},y_{2})$. Denoting $\hbox{div}_{x}w=\sum_{i=1}^{n}\partial_{x_{i}}w$, we derive
\begin{align*}\int_{\R^{N}}\hbox{div}_{x}[\partial_{y}v\nabla_{x}v]\, dx
=
\lim_{R\to+\infty}\int_{|x|\leq R}\hbox{div}_{x}[\partial_{y}v\nabla_{x}v]\, dx
=\lim_{R\to+\infty}\int_{|x|= R}\partial_{y}v\nabla_{x}v\cdot\tfrac x{|x|}\, d\sigma=0.\end{align*}
Therefore, multiplying (E) by $\partial_{y}v$ and integrating over $\R^{N}$ with respect to $x$, we obtain
\begin{align*}
0&=\int_{\R^{N}}-\partial^{2}_{y}v\partial_{y}v-\Delta_{x}v\partial_{y}v+v\partial_{y}v-f(v)\partial_{y}v\,dx\\
&=\int_{\R^{N}}-\tfrac 12\partial_{y}|\partial_{y}v|^{2}-\hbox{div}_{x}[\partial_{y}v\nabla_{x}v]+\tfrac 12
\partial_{y}|\nabla_{x}v|^{2}+\partial_{y}(\tfrac 12|v|^{2}-F(v))\,dx\\
&=\partial_{y}\,\left[\, \int_{\R^{N}}-\tfrac 12|\partial_{y}v|^{2}+
\tfrac 12|\nabla_{x}v|^{2}+\tfrac 12|v|^{2}-F(v)\,dx\,\right]\\
&=\partial_{y}\,\left[-\tfrac{1}{2}\|\partial_{y}v(\cdot,y)\|_{2}^{2}+\f(v(\cdot,y))\,\right]=-\partial_{y}E_{v}(y)
\end{align*}
and the Lemma follows.\QED
\subsection{The variational setting}
Fixed $b\in [0,c)$ we consider the space
\[ \XX_{b}=\{ v\in \MM\, /\, \liminf_{y\to
    \pm\infty}\hbox{dist}(v(\cdot,y),\VV^{b}_{\pm})= 
0\hbox{ and }
\inf_{y\in\R}\f(v(\cdot,y))\geq b\}
\]
on which we look for minima of the
functional
$$
\ff(v)=\int_{\R}
\tfrac{1}{2}\|\partial_{y}v(\cdot,y)\|_{2}^{2}+
	(\f(v(\cdot,y))-b) \,dy.
$$
\begin{Remark}\label{R:Mnonvuoto}{\rm
    The problem of finding a minimum of
    $\ff$ on $\XX_{b}$ is well posed. In fact, if $v\in\XX_{b}$ then $\f(v(\cdot,y))\geq b$
    for every $y\in\R$ and so the functional $\ff$ is well defined and non
    negative on $\XX_{b}$. Moreover $\XX_{b}\not=\emptyset$ and $m_{b}=\inf_{v\in\XX_{b}}\ff(v)<+\infty$.
    Indeed, for any $u\in \X$, recalling Lemma \ref{L:raggi} and considered
    the function    
    \[v(x,y)=\begin{cases}\omega_{b,u}u(x)& x\in\R^{N},\ y\geq\omega_{b,u},\\
    yu(x)& x\in\R^{N},\ \aub<y<\omega_{b,u},\\
    \aub u(x)& x\in\R^{N},\ y<\aub,\end{cases}\]
    we have that
    $v\in\XX_{b}$ and 
    $\ff(v)=\int_{\aub}^{\oub}\tfrac 12\|u\|^{2}_{2}+\f(yu)-b\,dy\leq
    (\tfrac 12\|u\|^{2}_{2}+\f(t_{u} u)-b)(\oub-\aub)<+\infty$.
}
\end{Remark}

\begin{Remark}\label{R:estensionef}{\rm 
    More generally, given an interval $I\subset\R$ we consider the functional
    \[
    \ff_{I}(v)=\int_{I}\tfrac{1}{2}\|\partial_{y}v(\cdot,y)\|_{2}^{2}+
	\f(v(\cdot,y))-b \,dy
    \]
    which is well defined for any $v\in \MM$
    such that $ \f(v(\cdot,y))\geq b$ for a.e. $y\in I$ or for every $v\in\MM$ if $I$ is bounded.}\end{Remark}

    We will make use of the following semicontinuity
    property
    \begin{Lemma}\label{L:semicontinuitafun}
    Let  $v\in \MM$ be such that $ \f(v(\cdot,y))\geq b$ for a.e. $y\in I\subset\R$. If $(v_{n})\subset\XX_{b}$ is such that $v_{n}\to v$ weakly in
    $H^{1}(S_{L})$ for any $L>0$, then
    $\ff_{I}(v)\leq\displaystyle{\liminf_{n\to\infty}}
    \ff_{I}(v_{n})$.\end{Lemma}
    \Proof Let $L_{1}<L_{2}\in\R$ be such that $(L_{1},L_{2})\subset I$. The sequence $(v_{n})$ is  weakly convergent to
 $v$ in   $H^{1}(S_{(L_1,L_2)})$ and constituted by radially symmetric functions in the $x$ variable. By Lemma III.2 in \cite{[L]} we derive that $v_{n}\to v$ strongly in $L^{p+1}(S_{(L_1,L_2)})$. Then, since by (\ref{eq:originf}) we have
$|F(v_{n})-F(v)|\leq \epsilon| v_{n}-v| (| v| +|v_{n}-v| )+2^{p-1}A_{\epsilon}
| v_{n}-v|( | v|^p+| v_{n}-v|^p)$,
we deduce that, as $n\to+\infty$,
\begin{align*}\int_{S_{(L_1,L_2)}} |F(v_{n})&-F(v)|\, dxdy
\leq\epsilon\|v_{n}-v\|_{L^{2}(S_{(L_1,L_2)})}(\|v\|_{L^{2}(S_{(L_1,L_2)})}+\|v_{n}-v\|_{L^{2}(S_{(L_1,L_2)})})\\
&
+
2^{p-1}A_{\epsilon}\|v_{n}-v\|_{L^{p+1}(S_{(L_1,L_2)})}(\|v\|_{L^{p+1}(S_{(L_1,L_2)})}^{p}+\|v_{n}-v\|_{L^{p+1}(S_{(L_1,L_2)})}^{p})\\
&\quad\quad=\epsilon\|v_{n}-v\|_{L^{2}(S_{(L_1,L_2)})}(\|v\|_{L^{2}(S_{(L_1,L_2)})}+\|v_{n}-v\|_{L^{2}(S_{(L_1,L_2)})})+o(1)
\end{align*}
Since $\epsilon$ is arbitrary and $(v_{n}-v)$ is bounded in $L^{2}(S_{(L_1,L_2)})$, we deduce $\int_{S_{(L_1,L_2)}} F(v_{n})\, dxdy\to\int_{S_{(L_1,L_2)}} F(v)dxdy$. Then, by the weak semicontinuity of the  norm, we obtain
\begin{align*}
\liminf_{n\to+\infty}\ff_{I}(v_{n})&\geq \liminf_{n\to+\infty}\ff_{(L_{1},L_{2})}(v_{n})
\\
&
=\liminf_{n\to+\infty}\tfrac 12\|v_{n}\|^{2}_{H^{1}(S_{(L_{1},L_{2})})}-\int_{S_{(L_{1},L_{2})}}F(v_{n})dx dy-b(L_{2}-L_{1})\\
&\geq \tfrac 12\|v\|^{2}_{H^{1}(S_{(L_{1},L_{2})})}-\int_{S_{(L_{1},L_{2})}}F(v)dxdy-b(L_{2}-L_{1})=
\ff_{(L_{1},L_{2})}(v)\end{align*}
and the Lemma follows by the arbitrariness of $L_{1}$ and $L_{2}$.\QED
\begin{Remark}\label{R:inequality}{\rm In the sequel we will study coerciveness properties of $\ff$. One of the key tools is the following simple estimate. Given $v\in \MM$, and $(y_{1},y_{2})\subset\R$ we have
     \begin{align*}
	    \ff_{(y_{1},y_{2})}(v)
	    &= \tfrac{1}{2} \int_{y_{1}}^{y_{2}}
	    \|\partial_{y}v(\cdot,y)\|_{2}^{2}\,dy
	    +\int_{y_{1}}^{y_{2}}\f(v(\cdot,y))-b\, dy\\  
	    &\geq \tfrac{1}{ 
	    2(y_{2}-y_{1})}\int_{\R^{N}}(\int_{y_{1}}^{y_{2}}
	    |\partial_{y}v(x,y)|\, dy)^{2}\, dx+\int_{y_{1}}^{y_{2}}\f(v(\cdot,y))-b\, dy\\  
	   &\geq \tfrac{1}{2(y_{2}-y_{1})}
	    \| 
	    v(\cdot,y_{1})-v(\cdot,y_{2})\|_{2}^{2}+\int_{y_{1}}^{y_{2}}\f(v(\cdot,y))-b\, dy.    \end{align*}
In particular if
    $\f(v(\cdot,y))\geq b+\nu$
    for any $y\in (y_{1},y_{2})$, then
    \begin{equation}\label{eq:stime2dim2} 	 \ff_{(y_{1},y_{2})}(v)
	   \geq \tfrac{1}{2(y_{2}-y_{1})}
	    \| 
	    v(\cdot,y_{1})-v(\cdot,y_{2})\|_{2}^{2}+\nu(y_{2}-y_{1})  \geq 
	    \sqrt{2\nu}\,\| v(\cdot,y_{1})-v(\cdot,y_{2})\|_{2}.
    \end{equation}
}\end{Remark}        

\begin{Remark}\label{R:stimabassam}{\rm In the sequel we will denote
$$
\delta_{0}=\delta((b+c)/2,(b+c)/2):=\hbox{dist}(\VV^{(b+c)/2}_{-},\VV^{(b+c)/2}_{+})\quad\hbox{and}\quad r_{0}=\frac{\delta_{0}}5
$$
By (\ref{eq:stime2dim2}) we can plainly prove that $m_{b}>0$. Indeed,
note that
if $v\in\XX_{b}$,  since by Lemma \ref{L:distanzaL2} we have $\delta_{0}>0$, by (\ref{eq:continuita1}), there exist $y_{1}<y_{2}\in\R$ such that
$\|v(\cdot,y_{1})-v(\cdot,y_{2})\|\geq \delta_{0}$ and $\f(v(\cdot,y))>(b+c)/2$ for any
$y\in (y_{1},y_{2})$. Then, by (\ref{eq:stime2dim2}) we obtain
$\ff_{(y_{1},y_{2})}(u)\geq
\sqrt{c-b}\,\delta_{0}>0$. In particular 
\[m_{b}\geq\sqrt{c-b}\,\delta_{0}.\]}\end{Remark}

One of the basic properties defining $\XX_{b}$ is the fact that if $v\in\XX_{b}$ then $\f(v(\cdot,y))\geq b$ for a.e. $y\in\R$. Unfortunately this condition is not necessarily preserved by the weak $H^{1}_{loc}$ convergence. To overcome this difficulty it is important the following Lemma.

\begin{Lemma}\label{L:st} Let $v\in\MM$ and $-\infty\leq \s<\tau\leq+\infty$ be such that
\begin{itemize}
\item[i)] $\f(v(\cdot,y))>b$ for any $y\in (\s,\tau)$
\item[ii)] either $\s=-\infty$ and $\liminf_{y\to-\infty}\dist(v(\cdot,y),\VV^{b}_{-})=0$ or
 $\s\in\R$ and $v(\cdot,\s)\in \VV^{b}_{-}$
\item[iii)] either $\tau=+\infty$ and $\liminf_{y\to+\infty}\dist(v(\cdot,y),\VV^{b}_{+})=0$ or
 $\tau\in\R$ and $v(\cdot,\tau)\in \VV^{b}_{+}$
\end{itemize}
then $\ff_{(\s,\tau)}(v)\geq m_{b}$. Moreover if $\liminf_{y\to \s^{+}}\f(v(\cdot,y))>b$ or $\liminf_{y\to \tau^{-}}\f(v(\cdot,y))>b$ then $\ff_{(\s,\tau)}(v)> m_{b}$.
\end{Lemma}
\Proof We consider the case in which $\s,\, \tau\in\R$. Similar arguments can be used to prove the statement in the cases $\s=-\infty$ or $\tau=+\infty$. 
We fix  two sequences $(s_{n}),\, (t_{n})\subset (\s,\tau)$ such that  $s_{n}\to \s$,
$t_{n}\to \tau$ as $n\to +\infty$ and
\begin{equation}\label{eq:ynpiumeno1} 
\f(v(\cdot,s_{n}))\leq\inf_{y\in (\s,s_{n})}\f(v(\cdot,y))+\tfrac{1}{n}\hbox{ and }\f(v(\cdot,t_{n}))\leq\inf_{y\in (t_{n},\tau)}\f(v(\cdot,y))+\tfrac{1}{n}.
\end{equation}
Moreover, since by (\ref{eq:continuita1}),we have $\|v(\cdot, s_{n})-v(\cdot,\s)\|_{2}\to 0$ and $\|v(\cdot, t_{n})-v(\cdot,\tau)\|_{2}\to 0$ as $n\to+\infty$, it is not restrictive to assume that
\begin{equation}\label{eq:ynpiumeno2}
\|v(\cdot, s_{n})-v(\cdot,\s)\|_{2}\leq r_{0}\hbox{ and }\|v(\cdot, t_{n})-v(\cdot,\tau)\|_{2}\leq r_{0}\hbox{ for any }n\in\N\end{equation}
For any $n\in\N$, consider the paths in $\X$ defined by 
\begin{align*}
 \gamma_{n,-}(y)&=v(\cdot,\s)+\tfrac{y-\s}{s_{n}-\s}\left(v(\cdot,s_{n})-v(\cdot,\s)\right),\quad y\in [\s,s_{n}],\\ 
\gamma_{n,+}(y)&=v(\cdot,\tau)+\tfrac{\tau-y}{\tau-t_{n}}\left(v(\cdot,t_{n})-v(\cdot,\tau)\right),\quad y\in [t_{n},\tau].
\end{align*}
Note that, for any $n\in\N$, the paths $\gamma_{n,-}$ and $\gamma_{n,+}$ continuously connect in $\X$ respectively the points $v(\cdot,\s)$, $v(\cdot,s_{n})$ and $v(\cdot,\tau)$, $v(\cdot,t_{n})$. Then, since by $(ii)-(iii)$, $\f(v(\cdot,\s)),\, \f(v(\cdot,\tau))\leq b$ and $\f(v(\cdot,s_{n})),\, \f(v(\cdot,t_{n}))> b$, defining for $n\in\N$
\begin{align*}
\bar s_{n}&=\inf\{\bar y\in [\s,s_{n}]\, /\, \f(\gamma_{n,-}(y))\geq b\hbox{ for any }y\in [\bar y, s_{n}]\},\\
\bar t_{n}&=\sup\{\bar y\in [t_{n},\tau]\, /\, \f(\gamma_{n,+}(y))\geq b\hbox{ for any }y\in [t_{n},\bar y]\},
\end{align*}
by continuity, we have that $\f(\gamma_{n,-}(\bar s_{n}))=b$ and 
$\f(\gamma_{n,+}(\bar t_{n}))=b$. Moreover, by definition, 
$\f(\gamma_{n,-}(y))\geq b$ for any $y\in [\bar s_{n},s_{n}]$ and
$\f(\gamma_{n,+}(y))\geq b$ for any $y\in [t_{n},\bar t_{n}]$.\\
Define, for $n,\, j\in\N$,
$$
w_{n,j}(\cdot,y)=\begin{cases}\gamma_{n,-}(\bar s_{n})&\hbox{ if }y\leq \bar s_{n},\\
\gamma_{n,-}(y)&\hbox{ if }\bar s_{n}<y\leq s_{n},\\
v(\cdot,y)&\hbox{ if }s_{n}<y\leq t_{j},\\
\gamma_{j,+}(y)&\hbox{ if }t_{j}<y\leq \bar t_{j},\\
\gamma_{j,+}(\bar t_{n})&\hbox{ if }\bar t_{j}<y,\end{cases}
$$
and note that $w_{n,j}\in\XX_{b}$, and so $\ff(w_{n,j})\geq m_{b}$ for any $n,\, j\in\N$.\\
To prove that $\ff_{(\s,\tau)}(v)\geq m_{b}$, 

we estimate the difference $\ff_{(\s,\tau)}(v)-\ff(w_{n,j})$. To this end, note that, since $\ff(w_{n,j})=\ff_{(\s,\tau)}(w_{n,j})$ and since
$w_{n,j}(\cdot,y)=v(\cdot,y)$ for any $y\in (s_{n},t_{j})$, we have
\begin{align}\nonumber
\ff_{(\s,\tau)}(v)&-\ff(w_{n,j})
=\int_{\s}^{s_{n}}
\tfrac 12\left(\|\partial_{y}v(\cdot,y)\|_{2}^{2}-\|\partial_{y}w_{n,j}(\cdot,y)\|_{2}^{2}\right)+
\left(\f(v(\cdot,y))-\f(w_{n,j}(\cdot,y))\right)\, dy\\ \label{eq:stimadalbasso1}&+\int_{t_{j}}^{\tau}
\tfrac 12\left(\|\partial_{y}v(\cdot,y)\|_{2}^{2}-\|\partial_{y}w_{n,j}(\cdot,y)\|_{2}^{2}\right)+
\left(\f(v(\cdot,y))-\f(w_{n,j}(\cdot,y))\right)\, dy.
\end{align}
Since $\partial_{y}w_{n,j}(\cdot,y)=\partial_{y}\gamma_{n,-}(\cdot,y)=\tfrac{1}{s_{n}-\s}\left(v(\cdot,s_{n})-v(\cdot,\s)\right)$ for $y\in (\bar s_{n},s_{n})$ and $\partial_{y}w_{n,j}(\cdot,y)=0$ for $y\in (\s,\bar s_{n})$, by (\ref{eq:continuita1}) we recover that
$$\int_{\s}^{s_{n}}
\|\partial_{y}w_{n,j}(\cdot,y)\|_{2}^{2}\, dy\leq \tfrac{1}{s_{n}-\s}\|v(\cdot,s_{n})-v(\cdot,\s)\|_{2}^{2}\leq\int_{\s}^{s_{n}}
\|\partial_{y}v(\cdot,y)\|_{2}^{2}\, dy.$$
Analogously, we obtain also that $\int_{t_{j}}^{\tau}\|\partial_{y}w_{n,j}(\cdot,y)\|_{2}^{2}dy\leq\int_{t_{j}}^{\tau}\|\partial_{y}v(\cdot,y)\|_{2}^{2}dy$ and by (\ref{eq:stimadalbasso1}) we conclude
\begin{equation}\label{eq:stimadalbasso2}
\ff_{(\s,\tau)}(v)-\ff(w_{n,j})\geq \int_{\s}^{s_{n}}\f(v(\cdot,y))-\f(w_{n,j}(\cdot,y))\, dy+\int_{t_{j}}^{\tau}\f(v(\cdot,y))-\f(w_{n,j}(\cdot,y))\, dy.
\end{equation}
Let us prove that $\liminf_{n\to +\infty} \int_{\s}^{s_{n}}\f(v(\cdot,y))-\f(w_{n,j}(\cdot,y))\, dy\geq 0$ for any $j\in\N$. Since $\f(w_{n,j}(\cdot,y))=\f(\gamma_{n,-}(y))$ for any $y\in (\bar s_{n},s_{n})$ and $\f(w_{n,j}(\cdot,y))=b$ for any $y\in (\s,\bar s_{n})$, we have
\begin{align*}
\int_{\s}^{s_{n}}&\f(v(\cdot,y))-\f(w_{n,j}(\cdot,y))\, dy=
\int_{\s}^{s_{n}}\f(v(\cdot,y))-\f(v(\cdot,s_{n}))\, dy\\&+
\int_{\s}^{\bar s_{n}}\f(v(\cdot,s_{n}))-b\, dy+
\int_{\bar s_{n}}^{s_{n}}\f(v(\cdot,s_{n}))-\f(\gamma_{n,-}(y))\, dy.\end{align*}
We separately estimate the three addenda at the right hand side of the above equality.\hb
For the first one that
by (\ref{eq:ynpiumeno1})
\begin{equation}\label{eq:addendo1}
\int_{\s}^{s_{n}}\f(v(\cdot,y))-\f(v(\cdot,s_{n}))\, dy\geq -\tfrac{1}{n}(s_{n}-\s).
\end{equation}
For the second one, we set 
$$
\ell:= \min\{\liminf_{y\to \s^{+}}\f(v(\cdot,y)),\, b+1\}
$$
and note that $\ell\in [b, b+1]$ and by (\ref{eq:ynpiumeno1}) we have in particular
$\ell\leq \liminf_{y\to \s^{+}}\f(v(\cdot,y))=\lim_{n\to+\infty}\f(v(\cdot,s_{n}))$
and hence that 
$\f(v(\cdot, s_{n}))-b\geq\ell-b+o(1)$ with $o(1)\to 0$ as $n\to +\infty$. Then
\begin{equation}\label{eq:addendo2}
\int_{\s}^{\bar s_{n}}\f(v(\cdot,y))-b\, dy\geq \left(\ell-b+o(1)\right)(\bar s_{n}- \s).
\end{equation}
Finally, for the third addendum, 
setting $4\mu=\ell-\f(v(\cdot,\s))$, 
we consider the two alternative cases: $\mu=0$ or $\mu>0$.\hb 
If $\mu=0$, since by $(ii)$, $\f(v(\cdot,\s))\leq b$ and $\ell\geq b$, we derive that $b=\ell=\f(v(\cdot,\s))$ and so, by (\ref{eq:ynpiumeno1}), that $\lim_{n\to+\infty}\f(v(\cdot,s_{n}))=\f(v(\cdot,\s))=b$. Then by Lemma \ref{L:taglio-conv-L2} and Remark
 \ref{R:conv-H1}  we derive $\bar s_{n}=\s$ and
$v(\cdot,s_{n})\to v(\cdot,\s)$ strongly in $\X$. Then, for any $y\in [\bar s_{n},s_{n}]$ we have $\|\gamma_{n,-}(y)-v(\cdot,\s)\|\leq\|v(\cdot,s_{n})-v(\cdot,\s)\|$ and by continuity of $\f$ we obtain $\sup_{y\in(\bar s_{n},s_{n})}\f(v(\cdot,s_{n}))-\f(\gamma_{n,-}(y))\to 0$ as $n\to+\infty$.  This allows us to conclude that if $\ell=\f(v(\cdot,\s))$ we have
\begin{equation}\label{eq:addendo3i}
\int_{\bar s_{n}}^{s_{n}}\f(v(\cdot,s_{n}))-\f(\gamma_{n,-}(y))\, dy\geq o(1)(s_{n}-\bar s_{n})\hbox{ as }n\to +\infty.
\end{equation}
In the second case, i.e. $\mu>0$, we have that necessarily $\liminf_{n\to+\infty}\|\nabla (v(\cdot,s_{n})-v(\cdot,s))\|\geq 8\mu_{0}$ for a certain $\mu_{0}>0$. Then
setting $\sigma_{n}(y)=\tfrac{y-\s}{s_{n}-\s}$ and $v_{n}=v(\cdot,s_{n})-v(\cdot,\s)$,  by Lemma \ref{L:taglio-conv-L2}, we obtain that for $n$ sufficiently large and $y\in (\bar s_{n},s_{n})$ we have
\begin{align*}\f(v(\cdot,s_{n}))-\f(\gamma_{n,-}(y))&=\f(v(\cdot,\s)+v_{n})-\f(v(\cdot,\s)+\sigma_{n}(y)v_{n})
\geq \mu_{0} (1-\sigma_{n}(y))=\mu_{0} \tfrac{s_{n}-y}{s_{n}-\s}.\end{align*}
Then
\begin{equation}\label{eq:addendo3ii}
\int_{\bar s_{n}}^{s_{n}}\f(v(\cdot,s_{n}))-\f(\gamma_{n,-}(y))\, dy\geq \mu_{0}\int_{\bar s_{n}}^{s_{n}}\tfrac{s_{n}-y}{s_{n}-\s}\, dy=\tfrac{\mu_{0}}{2}\tfrac{s_{n}-\bar s_{n}}{s_{n}-\s}(s_{n}-\bar s_{n}).\end{equation}
By (\ref{eq:addendo3i}) and (\ref{eq:addendo3ii}) we obtain
\begin{equation}\label{eq:addendo3finale}
\int_{\bar s_{n}}^{s_{n}}\f(v(\cdot,s_{n}))-\f(\gamma_{n,-}(y))\, dy\geq (\tfrac{\mu_{0}}{2}\tfrac{s_{n}-\bar s_{n}}{s_{n}-\s}+o(1))(s_{n}-\bar s_{n}).\end{equation}
Gathering (\ref{eq:addendo1}), (\ref{eq:addendo2}), (\ref{eq:addendo3finale}), we conclude that if $n$ is sufficiently large then, for any $j\in\N$,
\begin{align}\nonumber
 \int_{\s}^{s_{n}}\f(v(\cdot,y))-\f(w_{n,j}(\cdot,y))\, dy&\geq-\tfrac{1}{n}(s_{n}-\s)
+(\ell-c+o(1))(\bar s_{n}- \s)+\\ \label{eq:stima-}
&+(\tfrac{\mu_{0}}{2}\tfrac{s_{n}-\bar s_{n}}{s_{n}-\s}+o(1) )(s_{n}-\bar s_{n}).\end{align}
and then, $\liminf_{n\to +\infty} \int_{\s}^{s_{n}}\f(v(\cdot,y))-\f(w_{n,j}(\cdot,y))\, dy\geq 0$. 
In a symmetric way, we can prove that $\liminf_{j\to +\infty} \int_{t_{j}}^{\tau}\f(v(\cdot,y))-\f(w_{n,j}(\cdot,y))\, dy\geq 0$ for every $n\in\N$. 
Then, since $\ff(w_{n,j})\geq m_{b}$, by (\ref{eq:stimadalbasso2}) we conclude $\ff_{(\s,\tau)}(v)\geq m_{b}$.\\
Let us finally prove that if
$\liminf_{y\to \s^{+}}\f(v(\cdot,y))>b$
 then $\ff_{(\s,\tau)}(v)>m_{b}$.\\
Considering $\ell$ and $\mu_{0}$ defined as above we have $\ell>b$ and $\mu_{0}>0$ and hence
$2\tilde\mu:=\min\{\ell-b, \tfrac{\mu_{0}}{2}\}>0$. Since $
(\bar s_{n}- \s)^{2}+
(s_{n}-\bar s_{n})^{2}\geq\tfrac{1}{2}(s_{n}-s)^{2}$,
by (\ref{eq:stima-}) we obtain that for $n$  large and $j\in\N$
\begin{align*}
\int_{\s}^{s_{n}}\f(v(\cdot,y))-\f(w_{n,j}(\cdot,y))\, dy&\geq o(1)(s_{n}-\s)
+\tilde\mu\left[(\bar s_{n}- \s)+\tfrac{(s_{n}-\bar s_{n})^{2}}{s_{n}-\s}\right]\\ 
&\geq o(1) (s_{n}-\s)+\tilde\mu\left[\tfrac{(\bar s_{n}- \s)^{2}+(s_{n}-\bar s_{n})^{2}}{(s_{n}-\s)^{2}}\right](s_{n}-\s)\\
&\geq \tfrac{\tilde\mu}{4}(s_{n}-\s).
\end{align*}
Then, by (\ref{eq:stimadalbasso2}), since $\liminf_{j\to +\infty} \int_{t_{j}}^{\tau}\f(v(\cdot,y))-\f(w_{n,j}(\cdot,y))\, dy\geq 0$ for every $n\in\N$, we recover that for $n$ sufficiently large
$$
\ff_{(\s,\tau)}(v)-m_{b}\geq\liminf_{j\to+\infty}\,  [\ff_{(\s,\tau)}(v)-\ff(w_{n,j})]\geq \tfrac{\tilde\mu}{4}(s_{n}-s)>0.
$$
Simmetricaly we can prove that if $\liminf_{y\to \tau^{-}}\f(v(\cdot,y))>b$ then $\ff_{(\s,\tau)}(v)>m_{b}$.\QED

\subsection{ Estimates near the boundary of $\VV^{b}_{-}$ and $\VV_{+}^{b}$}
To study coercivity property of $\ff$ we first establish some technical local results.
We define the constants (depending on $b$)
\begin{equation}\label{eq:constantr0lambda0}
\ub=b+ \tfrac{c-b}4,\ \  
\hbox{ and }\ \Lambda_{0}=\sqrt{\tfrac{c-b}2}\, \tfrac {r_{0}}4
\end{equation}
where $\delta_{0}$ and $r_{0}$ are definined in Remark \ref{R:stimabassam}, noting that
\begin{equation}\label{eq:distanze}\mathrm{dist}(\VV^{b}_{-},\VV^{b}_{+})\geq \mathrm{dist}(\VV^{\ub}_{-},\VV^{\ub}_{+})\geq
5r_{0}.\end{equation}
Given $u_{0}\in \X$ we denote
\begin{align*}\XX^{-}_{b,u_{0}}&=\{v\in\MM\,/\,  v(\cdot,0)=u_{0},\ \inf_{(-\infty,0)}\f(v(\cdot,y))\geq b,\ 
\liminf_{y\to -\infty}\hbox{dist}(v(\cdot,y),\VV^{b}_{-})=0\},\\
\XX^{+}_{b,u_{0}}&=\{v\in\MM\,/\,  v(\cdot,0)=u_{0},\ \inf_{(0,+\infty)}\f(v(\cdot,y))\geq b,\ 
\liminf_{y\to +\infty}\hbox{dist}(v(\cdot,y),\VV^{b}_{+})=0\}.\end{align*}

Next Lemma establishes that if $v\in \XX^{+}_{b,u_{0}}$ (resp. $\XX^{-}_{b,u_{0}}$) is such that $\ff_{(0,+\infty)}(v)$ (resp. $\ff_{(-\infty, 0)}(v)$) is sufficiently small, then the trajectory $y\to v(\cdot,y)$ remains close to the set $\VV^{\ub}_{+}$ (resp. $\VV^{\ub}_{-}$) with respect to the $L^{2}(\R^{N})$ metric.
\begin{Lemma}\label{L:closel2}
If $u_{0}\in \X$, $\f(u_{0})\geq b$, $v\in\XX^{+}_{b,u_{0}}$ {\rm (}resp. $v\in\XX^{-}_{b,u_{0}}${\rm )} and
$\ff_{(0,+\infty)}(v)\leq\Lambda_{0}$ {\rm (}resp. $\ff_{(-\infty,0)}(v)\leq\Lambda_{0}${\rm )} then
$\mathrm{dist}(v(\cdot,y),\VV^{\ub}_{+})\leq r_{0}$  for every $y\in [0,+\infty)$ {\rm  (}resp. $\mathrm{dist}(v(\cdot,y),\VV^{\ub}_{+})\leq r_{0}$ for every $y\in (-\infty,0]${\rm )}.
\end{Lemma}
\Proof 
By (\ref{eq:continuita1}) the function $y\in [0,+\infty)\mapsto
v(\cdot,y)\in L^{2}(\R^{n})$ is continuous. Hence, using Remark \ref{R:V-bound+}, the map $y\in [0,+\infty)\mapsto\mathrm{dist}(v(\cdot,y),\VV^{\ub}_{+})$ is continuous too. 
If, by contradiction, $y_{0}\geq 0$ is such that $\mathrm{dist}(v(\cdot,y_{0}),\VV^{\ub}_{+})> r_{0}$, since $\liminf_{y\to+\infty}\mathrm{dist}(v(\cdot,y),\VV^{b}_{+})=0$, by continuity
there exists an interval $(y_{1},y_{2})\subset\R$ such that
$0<\mathrm{dist}(v(\cdot,y),\VV^{\ub}_{+})< r_{0}$ for any $y\in(y_{1},y_{2})$ and
$\|v(\cdot,y_{1})-v(\cdot,y_{2})\|_{2}\geq r_{0}/2$. By (\ref{eq:distanze}) we derive
$v(\cdot,y)\notin \VV^{\ub}_{+}\cup \VV^{\ub}_{-}$ and so
$\f(v(\cdot,y))-b\geq \ub-b=(c-b)/4$ for all $y\in (y_{1},y_{2})$ . By (\ref{eq:stime2dim2}) we conclude
\[\Lambda_{0}\geq\ff_{(0,+\infty)}(v)\geq\ff_{(y_{1},y_{2})}(v)\geq \sqrt{\tfrac{c-b}2}\|v(\cdot,y_{1})-v(\cdot,y_{2})\|_{2}\geq\sqrt{\tfrac{c-b}2}\, \tfrac{r_{0}}2=2\Lambda_{0},\]
a contradiction which proves the Lemma. Analogous is the proof in the case $v\in\XX^{-}_{b,u_{0}}$.\QED

Clearly the infimum value of $\ff_{(0,+\infty)}$ on $\XX^{+}_{b,u_{0}}$ is close to $0$ as $\dist(u_{0},\VV^{b})$ is small. Next result displays  a test function $w^{+}_{u_{0}}\in\XX^{+}_{b,u_{0}}$ which gives us more precise information.
\begin{Lemma}\label{L:minimonormalim}
Let $b\in[0,c)$, then there exists $C_{+}(b)>0$ such that for every $u_{0}\in\VV^{\ub}_{+}\setminus \VV^{b}_{+}$ there exists $w^{+}_{u_{0}}\in\XX^{+}_{b,u_{0}}$ such that
\begin{align*}\sup_{y>0}\|w^{+}_{u_{0}}(\cdot,y)-u_{0}\|_{2}\leq  \tfrac1 {\nu^+(\ub)}(\f(u_{0})-b)
\ \hbox{ 
and}\quad
\ff_{(0,+\infty)}(w^{+}_{u_{0}})\leq C_{+}(b)(V(u_{0})-b)^{3/2}.\end{align*}
\end{Lemma}
\Proof Note that, since $u_{0}\in\VV^{\ub}_{+}$, by Lemma \ref{L:raggi}, we have $\f'(u_{0})u_{0}<0$ and there exists a unique $s_{0}\in (1,+\infty)$ such that $V(su_{0})>b$ for any $s\in [1,s_{0})$ and $V(s_{0}u_{0})=b$. Moreover
$\tfrac{d}{ds}\f(su_{0})=s(\f'(u_{0})u_{0}+\int_{\R^{N}}f(u_{0})u_{0}-\tfrac 1sf(su_{0})u_{0}\, dx)$
and since, by (\ref{eq:convessita}),  $\int_{\R^{N}}f(u_{0})u_{0}-\tfrac 1sf(su_{0})u_{0}\, dx\leq0$ for any $s\geq1$, we deduce that $\tfrac{d}{ds}\f(su_{0})\leq s\f'(u_{0})u_{0}$ for any $s\geq 1$.
Integrating this last inequality on the interval $[1,s_{0}]$, we obtain
$\f(s_{0}u_{0})\leq\f(u_{0})+\tfrac12(s^{2}_{0}-1)\f'(u_{0})u_{0}$ and so the estimate $s_{0}-1\leq \tfrac{\f(u_{0})-b}{|\f'(u_{0})u_{0}|}$.  We define
\[w^{+}_{u_{0}}(x,y)=\begin{cases}u_{0}(x)&y\leq 0\\ 
(1+\tfrac{y^{2}}2)u_{0}&y\in (0,{\scriptstyle\sqrt{2(s_{0}-1)}})\\
s_{0}u_{0}&y\geq {\scriptstyle\sqrt{2(s_{0}-1)}}\end{cases}\]
notiing that $w^{+}_{u_{0}}\in\XX^{+}_{b,u_{0}}$ and 
$\sup_{y\geq 0}\|w^{+}_{u_{0}}(\cdot,y)-u_{0}\|_{2}=(s_{0}-1)\|u_{0}\|_{2}\leq \tfrac{\f(u_{0})-b}{|\f'(u_{0})u_{0}|}\|u_{0}\|_{2}$. 
Moreover, since $s_{0}-1\leq \tfrac {\f(u_{0})-b}{|\f'(u_{0})u_{0}|}$, we get
\begin{align*}\ff_{(-\infty,0)}(w^{+}_{u_{0}})&=\int_{0}^{{\scriptstyle\sqrt{2(s_{0}-1)}}}\tfrac 12\|\partial_{y}(1+\tfrac {y^{2}}2)u_{0}(\cdot)\|_{2}^{2}\, dy+\int_{0}^{{\scriptstyle\sqrt{2(s_{0}-1)}}}\f((1+\tfrac {y^{2}}2)u_{0}(\cdot))-b\, dy\\
&\leq\int_{0}^{{\scriptstyle\sqrt{2(s_{0}-1)}}}\tfrac 12y^{2}\|u_{0}\|_{2}^{2}\, dy+\int_{0}^{{\scriptstyle\sqrt{2(s_{0}-1)}}}\f(u_{0})-b\, dy\\
&= {\scriptstyle\sqrt{2(s_{0}-1)}}(\tfrac {(s_{0}-1)}3\|u_{0}\|_{2}^{2}+(\f(u_{0})-b))\\
&\leq{\scriptstyle\sqrt{\tfrac{2}{|\f'(u_{0})u_{0}|}}}(\tfrac1 {3|\f'(u_{0})u_{0}|}\|u_{0}\|_{2}^{2}+1)(\f(u_{0})-b)^{3/2}.
\end{align*}
By Lemma \ref{L:gradienteV+} we know that
$|\f'(u_{0})u_{0}|\geq \nu^{+}(\ub)\max\{1,\|u_{0}\|^{2}_{2}\}$  and the Lemma follows  considering
$C_{+}(b)= {\scriptstyle\sqrt{\tfrac{2}{\nu^{+}(\ub)}}}(\tfrac1 {3\nu^+(\ub)}+1)$. \QED

\noindent For any $b\in [0,c)$ we fix $\bi_{+}\in(b,\ub]$ such that the following inequalities hold true:
\begin{equation}\label{eq:b*}
\tfrac{\bi_{+}-b}{\nu^{+}(\ub)}<\tfrac 12,\quad \max\{1,C_{+}(b)\}(\bi_{+}-b)^{1/4}<\tfrac 14,\quad C_{+}(b)(\bi_{+}-b)^{3/2}\leq\Lambda_{0}.\end{equation}
Next Gronwall type result  will play an important role together with Lemma \ref{L:minimonormalim}. 
\begin{Lemma}\label{L:finitetime}
Assume that  $u_{0}\in\VV^{\bi_{+}}_{+}\setminus\VV^{ b}_{+}$ and $v\in\XX_{b,u_{0}}^{+}$ are such that  
\begin{equation}\label{eq:coerc1}\hbox{ if }y\in [0,1)\hbox{ is such that }V(\bar v(\cdot,y))\leq \bi_{+}\hbox{ then }\ff_{(y,+\infty)}(\bar v)\leq C_{+}(b)(\f(\bar v(\cdot, y))-b)^{3/2}.\end{equation}
Then  there exists $\bar y\in (0,1)$ such that $V(v(\cdot, \bar y))=b$, $v(\cdot, \bar y)\in\VV^{b}_{+}$ and
$v(\cdot,y)=v(\cdot,\bar y)$ for every $y\in[\bar y,+\infty).$
\end{Lemma}
\Proof We first note that, since $u_{0}\in\VV^{\bi_{+}}_{+}\setminus\VV^{ b}_{+}$ and $v\in\XX_{b,u_{0}}^{+}$ we have $V(v(\cdot,0))=V(u_{0})\le \bi_{+}$ and hence, by (\ref{eq:coerc1}) and (\ref{eq:b*}), we have
$\ff_{(0,+\infty)}(v)\leq C_{+}(b)(V(u_{0})-b)^{3/2}\leq\Lambda_{0}$.  By Lemma \ref{L:closel2} we then deduce that $\mathrm{dist}(v(\cdot,y),\VV^{\beta}_{+})\leq r_{0}$ for any $y>0$ and, by the definition of $r_{0}$, we obtain that $v(\cdot,y)\notin \VV^{\bi_{+}}_{-}$ for any $y>0$. In particular, if $y>0$ and $V(v(\cdot,y))\leq \bi_{+}$ then $v(\cdot,y)\in\VV^{\bi_{+}}_{+}$.\\ We claim that there exists a sequence $(\zeta_{n})\subset [0,\frac12)$ such that 
\begin{equation}\label{eq:induct}\zeta_{n-1}<\zeta_{n}\leq\zeta_{n-1}+(\tfrac {\bi_{+}-b}{4^{2(n-1)}})^{1/4}<\tfrac12\hbox{ and
}\f(v(\cdot,\zeta_{n}))-b\leq\tfrac {\bi_{+}-b}{4^{n}}, \quad\forall n\in\N.\end{equation}
Indeed, defining
$\zeta_{0}=0$
 by (\ref{eq:b*}) and (\ref{eq:coerc1}) we have that for any $\zeta>\zeta_{0}$ 
\begin{align*}
\int_{\zeta_{0}}^{\zeta}\f(v(\cdot,s))-b\, ds&\le \ff_{(\zeta_{0},+\infty)}(v)\leq C_{+}(b)(\f(v(\cdot,\zeta_{0}))-b)^{3/2}\\
&\leq
C_{+}(b)(\bi_{+}-b)^{1/4}(\bi_{+}-b)(\bi_{+}-b)^{1/4}\leq \tfrac 1{4}(\bi_{+}-b)(\bi_{+}-b)^{1/4},
\end{align*} 
and so 
\begin{equation}\label{eq:din}
\exists\, \zeta_{1}\in (\zeta_{0},\zeta_{0}+(\bi_{+}-b)^{1/4})\hbox{ such that } \f(\bar v(\cdot,\zeta_{1}))-b\leq\tfrac{\bi_{+}-b}{4},
\end{equation}
Note that,  by (\ref{eq:b*}), $\zeta_{0}+(\bi_{+}-b)^{1/4}<\zeta_{0}+\tfrac1{4}<\tfrac12$ and so $\zeta_{1}\in (0,\frac12)$.\\
Now, if $\zeta_{n}$  verifies (\ref{eq:induct}) by
 (\ref{eq:coerc1}) we obtain  that for any
$\zeta>\zeta_{n}$ 
\begin{align*}\int_{\zeta_{n}}^{\zeta}\f( v(\cdot,s))-b\, ds&\le \ff_{(\zeta_{n},+\infty)}(v)\leq C_{+}(b)(\f(v(\cdot,\zeta_{n}))-b)^{3/2}\\
&\leq
C_{+}(b)(\bi_{+}-b)^{1/4}(\tfrac{\bi_{+}-b}{4^{n}})(\tfrac{\bi_{+}-b}{4^{2n}})^{1/4}< \tfrac{\bi_{+}-b}{4^{n+1}}(\tfrac{\bi_{+}-b}{4^{2n}})^{1/4},\end{align*}
implying that
$$
\exists\, \zeta_{n+1}\in (\zeta_{n},\zeta_{n}+(\tfrac{\bi_{+}-b}{4^{2n}})^{1/4})\hbox{ such that } \f(v(\cdot,\zeta_{n+1}))-b\leq\tfrac {\bi_{+}-b}{4^{n+1}},
$$
with, by (\ref{eq:b*}),
$$
\zeta_{n+1}<\sum_{j=0}^{n}(\tfrac{\bi_{+}-b}{4^{2j}})^{1/4}=(\bi_{+}-b)^{1/4}\sum_{j=0}^{+\infty}\tfrac 1{2^{j}}<\tfrac12.
$$
Then, by induction, (\ref{eq:induct}) holds true for any $n\in\N$.\\
 Now,
note that  by (\ref{eq:induct}) we have $\zeta_{n}\to \bar y\in (0,\frac12]$ as $n\to +\infty$. Moreover, since $v\in\XX_{b,u_{0}}$ there result $\f( v(\cdot,\zeta_{n}))\geq b$ for all $n\in\N$ and hence, by (\ref{eq:induct}), $\f( v(\cdot,\zeta_{n}))\to b$. Then, by Lemma \ref{L:semicVv}, we deduce  $\f(v(\cdot,\bar y))=b$.
Moreover, by (\ref{eq:continuita1}), $v(\cdot,\zeta_{n})\to v(\cdot,\bar y)$ in $L^{2}(\R^{N})$. Then we can conclude that  $v(\cdot,\bar y)\in\VV^{b}_{+}$ and hence, using (\ref{eq:coerc1}), that $\ff_{(\bar y,+\infty)}(v)\leq C_{+}(b)(\f(v(\cdot,\bar y))-b)^{3/2}=0$, which implies $v(\cdot,y)= v(\cdot,\bar y)$ for every $y\geq\bar y$.
\QED

Lemma \ref{L:finitetime} and Lemma \ref{L:minimonormalim} have in particular the  following consequence which will be a key tool in constructing minimizing sequences for $\ff$ with suitable compactness properties.
\begin{Lemma}\label{L:dinamicaparallela} Let $b\in [0,c)$ then, for every $u_{0}\in\VV^{\bi_{+}}_{+}\setminus\VV^{ b}_{+}$ and $v\in\XX^{+}_{b,u_{0}}$
there exists $\tilde v\in \XX^{+}_{b,u_{0}}$ such that
\[\sup_{y\in (0,+\infty)}\|\tilde v(\cdot,y)-u_{0}\|_{2}\leq 1\hbox{ and }\ff_{(0,+\infty)}(\tilde v)\leq\min\{\Lambda_{0},\ff_{(0,+\infty)}(v)\}.\]
\end{Lemma}
\Proof 
Note that, by Lemma \ref{L:minimonormalim} and (\ref{eq:b*}), we have in that if $u_{0}\in\VV^{\bi_{+}}_{+}\setminus \VV^{b}_{+}$ then $\ff_{(0,+\infty)}(w^{+}_{u_{0}})\leq\Lambda_{0}$ and $\|w^{+}_{u_{0}}(\cdot,y)-u_{0}\|\leq\tfrac 12$ for any $y>0$. In particular if $u_{0}\in\VV^{\bi_{+}}_{+}\setminus \VV^{b}_{+}$
and $v\in\XX^{+}_{b,u_{0}}$ are such that $\ff_{(0,+\infty)}(v)>\Lambda_{0}$ then the statement of the Lemma holds true with $\tilde v=w^{+}_{u_{0}}$.\par\noindent 
To prove the Lemma we argue by contradiction assuming that there exist $u_{0}\in\VV^{\bi_{+}}_{+}\setminus\VV^{ b}_{+}$ and 
$v\in \XX^{+}_{b,u_{0}}$ with $\ff_{(0,+\infty)}(v)\leq\Lambda_{0}$ such that 
\begin{equation}\label{eq:contradiction}
\ff_{(0,+\infty)}(\tilde v)>\ff_{(0,+\infty)}( v)
\hbox{ for every }
\tilde v\in \XX^{+}_{b,u_{0}}\hbox{ such that }\sup_{y\in (0,+\infty)}\|\tilde v(\cdot,y)-u_{0}\|_{2}\leq 1.
\end{equation}  
By (\ref{eq:contradiction}) we have $\sup_{y\in (0,+\infty)}\|v(\cdot,y)-u_{0}\|_{2}> 1$ and since $v(\cdot,0)=u_{0}$, by (\ref{eq:continuita1}) we recover that
\begin{equation}\label{eq:y0}
\exists\, y_{0}>0\hbox{ such that }\| v(\cdot,y_{0})-u_{0}\|_{2}=\tfrac 12\hbox{ and }\|v(\cdot,y)-u_{0}\|_{2}<\tfrac 12\hbox { for any }y\in [0,y_{0}).
\end{equation} 
As already noted in the proof of the previous Lemma, by Lemma \ref{L:closel2}, since
$\ff_{(0,+\infty)}(v)\leq \Lambda_{0}$,   we have that if $y>0$ and $V( v(\cdot,y))\leq \bi_{+}$ then $v(\cdot,y)\in\VV^{\bi_{+}}_{+}$.
We deduce that
\begin{equation}\label{eq:coerc}
\hbox{ if }\tilde y\in [0,y_{0})\hbox{ and }V(v(\cdot,\tilde y))\leq \bi_{+}\hbox{ then }\ff_{(\tilde y,+\infty)}( v)\leq C_{+}(b)(\f(v(\cdot, \tilde y))-b)^{3/2}.
\end{equation}
Indeed, considering the function
$$
 \tilde v(\cdot,y)=\begin{cases} v(\cdot,y)&0\leq y< \tilde y\\
w^{+}_{ v(\cdot,\tilde y)}(\cdot,y-\tilde y)&y\geq \tilde y,\end{cases}
$$
we have $\tilde v\in\XX^{+}_{b,u_{0}}$. Now note that for every $y\in[0,\tilde y)\subset[0,y_{0})$, by definition of $y_{0}$ we have
$
\|\tilde v(\cdot,y)-u_{0}\|_{2}=\| v(\cdot,y)-u_{0}\|_{2}<\frac12
$
while if $y\geq \tilde y$, then by Lemma \ref{L:minimonormalim} and (\ref{eq:b*})
$$
\|\tilde v(\cdot,y)-u_{0}\|_{2}=\|w^{+}_{ v(\cdot,\tilde y)}(\cdot,y-\tilde y)-u_{0}\|_{2}\le\frac{\bi_{+}-b}{\nu^{+}(\bi_{+})}<\frac12.
$$
Hence we recover that
$\sup_{y>0}\|\tilde v(\cdot,y)-u_{0}\|_{2}\le 1$.
Then, by (\ref{eq:contradiction})  we obtain $\ff_{(0,+\infty)}( v)<\ff_{(0,+\infty)}(\tilde v)\le\ff_{(\tilde y,+\infty)}(\tilde v)=\ff_{(0,+\infty)}(w^{+}_{v(\cdot,\tilde y)})$ and  (\ref{eq:coerc}) follows by Lemma \ref{L:minimonormalim}.\\
We now note that, by Remark \ref{R:inequality} we have $\ff_{(0,y_{0})}(v)\geq\tfrac 1{2y_{0}}\|v(\cdot,y_{0})-u_{0}\|_{2}^{2}= \tfrac 1{8y_{0}}$ and so, by (\ref{eq:b*}) and (\ref{eq:coerc}), we deduce
$y_{0}\geq \tfrac 1{8C_{+}(b)(\bi_{+}-b)^{3/2}}>1$. Then, by (\ref{eq:coerc}) and Lemma \ref{L:finitetime}, we derive that there exists
$\bar y\in (0,1)$ such that $v(\cdot,\bar y)\in\VV^{b}_{+}$ and 
$v(\cdot, y)= v(\cdot,\bar y)$ for any $y\geq \bar y$.
Hence, using  (\ref{eq:y0}), we obtain $1<\sup_{y\in (0,+\infty)}\| v(\cdot,y)-u_{0}\|_{2}=\sup_{y\in (0,\bar y]}\|\ v(\cdot,y)-u_{0}\|_{2}\le 
\sup_{y\in (0,y_{0}]}\| v(\cdot,y)-u_{0}\|_{2}=\tfrac 12$,
a contradiction which proves the Lemma.\QED

\noindent The following Lemma is an analogous of Lemma \ref{L:minimonormalim} for $\XX^{-}_{b,u_{0}}$ when $b>0$. We omit the proof since it is bases on an argument  symmetric to the one used proving Lemma \ref{L:minimonormalim}, using  Lemma \ref{L:gradienteV-} instead of Lemma \ref{L:gradienteV+}.

\begin{Lemma}\label{L:minimonormalim-}
Let $b\in(0,c)$, then there exists $C_{-}(b)>0$ such that for any $u_{0}\in\VV^{\ub}_{-}\setminus \VV^{b}_{-}$ there exists $ w^{-}_{u_{0}}\in\XX^{-}_{b,u_{0}}$ such that
$$
\ff_{(-\infty,0)}(w^{-}_{u_{0}})\leq C_{-}(b)(V(u_{0})-b)^{3/2}.
$$
\end{Lemma}

\noindent For any $b\in (0,c)$ we fix $\bi_{-}\in(b,\ub]$ such that the following inequalities hold true:
\begin{equation}\label{eq:b*-}
\max\{1,C_{-}(b)\}(\bi_{-}-b)^{1/4}<\tfrac 14\quad\hbox{and}\quad C_{-}(b)(\bi_{-}-b)^{3/2}\leq\Lambda_{0}.\end{equation}
Analogously to Lemma \ref{L:finitetime} we can prove
\begin{Lemma}\label{L:finitetime-}
Let $b\in(0,c)$ and assume that $u_{0}\in\VV^{\bi_{-}}_{-}\setminus\VV^{ b}_{-}$ and $v\in\XX_{b,u_{0}}^{-}$ are such that  
\begin{equation}\label{eq:coerc1-}\hbox{ if }y\in (-1,0]\hbox{ is such that }V( v(\cdot,y))\leq \bi_{-}\hbox{ then }\ff_{(-\infty,y)}( v)\leq C_{-}(b)(\f(v(\cdot, y))-b)^{3/2}.\end{equation}
Then, there exists $\bar y\in (-1,0)$ such that $V(v(\cdot, \bar y))=b$, $v(\cdot, \bar y)\in\VV^{b}_{-}$ and
$v(\cdot,y)= v(\cdot,\bar y)$ for any $y\in(-\infty,\bar y].$
\end{Lemma}

\noindent The situation is slightly different when $b=0$.
\begin{Lemma}\label{L:minimonormalim0-}
If $b=0$ there exists $\bi_{0}\in (0,\frac c4)$ such that for any
$u_{0}\in\VV^{\bi_{0}}_{-}\setminus \{0\}$ there exists  $w^{-}_{u_{0}}\in\XX^{-}_{b,u_{0}}$ such that $\ff_{(-\infty,0)}(w^{-}_{u_{0}})\leq 3\f(u_{0})$.
\end{Lemma}
\Proof If $u_{0}\in\VV^{\bi_{0}}_{-}$ for some $\bi_{0}\in (0,\frac c4)$ we set
$$
w^{-}_{u_{0}}(x,y)=\begin{cases}u_{0}(x)&y\geq 0\\ 
(1+y)u_{0}(x)&y\in (-1,0)\\
0&y\leq -1\end{cases}
$$
noting that $w^{-}_{u_{0}}\in \XX^{-}_{0,u_{0}}$ and
$\ff_{(-\infty,0)}(w^{-}_{u_{0}})\leq
\int_{-1}^{0}\tfrac 12\|u_{0}\|_{2}^{2}+V(u_{0})\, dy\leq \tfrac 12\|u_{0}\|_{2}^{2}+V(u_{0}).$
By Remark \ref{R:V-bound} and Lemma \ref{L:origine}, if $\bi_{0}$ is sufficiently small,  we obtain $ \|u_{0}\|_{2}^{2}\le 4V(u_{0})$ and  so
$\ff_{(-\infty,0)}(w^{-}_{u_{0}})\leq 3V(u_{0})$.\QED
\begin{Remark}\label{R:min-lambda0} {\rm Eventually taking $\bi_{0}$ smaller, we can assume that
$\ff_{(-\infty,0)}(w^{-}_{u_{0}})\leq\Lambda_{0}$ for  $u_{0}\in\VV^{\bi_{0}}_{-}$}
\end{Remark}

\subsection{Minimizing sequences and their limit points}
The local results that we have described in the previous section, allow us to produce a minimizing sequence of $\ff$ on $\XX_{b}$ with suitable compactness properties.
\begin{Lemma}\label{L:minseq}  For every $b\in [0,c)$ there exists $L_{0}>0$, $\bar C>0$ and $(v_{n})\subset\XX_{b}$ such that $\ff(v_{n})\to m_{b}$ and
\begin{description}
\item{i)} $\mathrm{dist}(v_{n}(\cdot,y),\VV^{\beta}_{-})\leq r_{0}$ for any $y\leq 0$ and $n\in\N$,
\item{ii)} $\mathrm{dist}(v_{n}(\cdot,y),\VV^{\beta}_{+})\leq r_{0}$ for any $y\geq L_{0}$ and $n\in\N$,
\item{iii)} $\|v_{n}(\cdot,y)\|_{2}\leq \bar C$ for any $y\in\R$ and $n\in\N$
\item{iv)} for every bounded interval $(y_{1},y_{2})\subset\R$ there exists $\hat C>0$, depending only on $y_{2}-y_{1}$, such that $\|v_{n}\|_{H^{1}(S_{(y_{1},y_{2})})}\leq \hat C$.
\end{description}
\end{Lemma}
\Proof Let $b\in [0,c)$ and $(w_{n})\subset\XX_{b}$ be such that $\ff(w_{n})\leq m_{b}+1$ for any $n\in\N$ and $\ff(w_{n})\to m_{b}$. We denote $\bi^{*}=\min\{\bi_{-},\bi_{+}\}$. Let 
$s_{n}=\sup\{y\in\R\, |\, \ff_{(-\infty,y)}(w_{n})\leq\Lambda_{0}\}$ and note that by Remark \ref{R:stimabassam}, since $\Lambda_{0}<m_{b}\leq\ff(w_{n})$, we have $s_{n}\in\R $ and
$\ff_{(-\infty,s_{n})}(w_{n})=\Lambda_{0}$. Since $w_{n}(\cdot,\cdot+s_{n})\in\XX^{-}_{b,w_{n}(\cdot,s_{n})}$ and $\ff_{(-\infty,0)}(w_{n}(\cdot,\cdot+s_{n}))=\Lambda_{0}$, by Lemma \ref{L:closel2} we derive that
$\mathrm{dist}(w_{n}(\cdot,y+s_{n}),\VV^{\beta}_{-})\leq r_{0}$ for any $y\leq 0$ and so, by (\ref{eq:distanze}), $\mathrm{dist}(w_{n}(\cdot,y),\VV^{b^{*}}_{+})\geq 4r_{0}$ for any $y\leq s_{n}$. We conclude that if $y\leq s_{n}$ and $\f(w_{n}(\cdot,y))\leq b^{*}$ then $w_{n}(\cdot,y)\in\VV^{b^{*}}_{-}$. 
A symmetric argument shows that there exists $t_{n}>s_{n}$ such that if
$y\geq t_{n}$ and $\f(w_{n}(\cdot,y))\leq b^{*}$ then $w_{n}(\cdot,y)\in\VV^{b^{*}}_{+}$. Define now 
\[y^{-}_{n}=\sup\{y\leq t_{n}\,|\, w_{n}(\cdot,y)\in \VV^{b^{*}}_{-}\}\hbox{ and }y_{n}^{+}=\inf\{ y\geq y_{n}^{-}\,|\, w_{n}(\cdot,y)\in \VV^{b^{*}}_{+}\}.\]  
Since $\liminf_{y\to\pm\infty}\f(w_{n}(\cdot,y))=b>\bi^{*}$ we deduce that $y^{-}_{n},y^{+}_{n}\in\R$.\\
Using Remarks  \ref{R:V-bound} and  \ref{R:V-bound+} we also recognize that $w_{n}(\cdot,y^{-}_{n})\in\VV^{b^{*}}_{-}$ and $w_{n}(\cdot,y^{+}_{n})\in\VV^{b^{*}}_{+}$.
Since the function $y\mapsto w_{n}(\cdot,y)$ is continuous with respect to the $L^{2}(\R^{N})$ metric and
$\mathrm{dist}(\VV^{\bi^{*}}_{-},\VV^{\bi^{*}}_{+})\geq 5r_{0}$ we deduce $y_{n}^{-}<y_{n}^{+}$.
Moreover $V(w_{n}(\cdot,y))>\bi^{*}$ for any $y\in (y_{n}^{-},y_{n}^{+})$ and $\|w_{n}(\cdot,y_{n}^{+})-w_{n}(\cdot,y_{n}^{-})\|_{2}\geq 5r_{0}$. By (\ref{eq:stime2dim2}) we derive
\begin{equation}\label{eq:y+-} y_{n}^{+}-y_{n}^{-}\leq \tfrac{\ff_{(y_{n}^{-},y_{n}^{+})}(w_{n})}{\bi^{*}-b}\leq\tfrac{m_{b}+1}{\bi^{*}-b}:= L_{0}\hbox{ and }\sup_{y\in (y_{n}^{-},y_{n}^{+}]}\|w_{n}(\cdot,y)-w_{n}(\cdot,y_{n}^{-})\|_{2}\leq \tfrac{m_{b}+1}{\scriptstyle{\sqrt{2(\bi^{*}-b)}}}.\end{equation}
We now claim that, eventually modifying the function $w_{n}$ on the set $\R^{N}\times [ (-\infty,y_{n}^{-})\cup (y_{n}^{+},+\infty)]$,  $w_{n}$ satisfies \begin{description}
\item{(I)} $\ff_{(-\infty,y_{n}^{-})}(w_{n})\leq\Lambda_{0}$,
\item{(II)} $\ff_{(y_{n}^{+},+\infty)}(w_{n})\leq\Lambda_{0}\hbox{ and }\|w_{n}(x,y)-w_{n}(x,y_{n}^{+})\|_{2}\leq 1$ for any $y\geq y_{n}^{+}$.
\end{description}
Indeed, 
 if (I) is not satisfied, since $w_{n}(\cdot,y_{n}^{-})\in\VV^{\bi_{-}}_{-}$, we 
can consider the new function
\[ w^{*}_{n}(\cdot,y)=\begin{cases}w^{-}_{w_{n}(\cdot,y_{n}^{-})}(\cdot,y-y_{n}^{-})&\hbox{if } y\leq y_{n}^{-}\\
w_{n}(\cdot,y)&\hbox{if }y>y_{n}^{-}.\end{cases}\]
nothing that $w^{*}_{n}\in \XX_{b}$, $\ff(w^{*}_{n})\leq\ff(w_{n})$ and $w^{*}_{n}$ satisfies (I)  by Lemma \ref{L:minimonormalim-}, Lemma \ref{L:minimonormalim0-}, (\ref{eq:b*-}) and Remark \ref{R:min-lambda0} .\par\noindent
Now,  assuming that (I) is verified, if (II) is not satisfied, since $w_{n}(\cdot,y_{n}^{+})\in\VV^{\bi_{+}}_{+}$ and $w_{n}(\cdot,\cdot+y_{n}^{+})\in\XX^{+}_{b,w_{n}(\cdot,y_{n}^{+})}$, by Lemma
\ref{L:dinamicaparallela} there exists a function $\tilde w_{n}\in\XX^{+}_{b,w_{n}(\cdot,y_{n}^{+})}$ such that $\ff_{(y_{n}^{+},+\infty)}(\tilde w_{n}(\cdot,\cdot-y_{n}^{+}))\leq\min\{\Lambda_{0},
\ff_{(y_{n}^{+},+\infty)}( w_{n})\}$ and $\|\tilde w_{n}(\cdot,y-y_{n}^{+})-w_{n}(x,y_{n}^{+})\|_{2}\leq 1$ for any $y\geq y_{n}^{+}$. Then considering
$$
w^{**}_{n}(\cdot,y)=\begin{cases}\tilde w_{n}(\cdot,y-y_{n}^{+})&\hbox{if } y\geq y_{n}^{+}\\
w_{n}(\cdot,y)&\hbox{if }y<y_{n}^{+}.\end{cases}
$$
 we recognize that $w^{**}_{n}\in \XX_{b}$, $\ff(w^{**}_{n})\leq\ff(w_{n})$ and $w^{**}_{n}$ satisfies (I) and (II). Hence, eventually modifying $w_{n}$ as indicated above our claim follows.\par\noindent
We finally set $v_{n}=w_{n}(\cdot,\cdot+y_{n}^{-})$ obtaing that $v_{n}\in\XX_{b}$ and $\ff(v_{n})=\ff(w_{n})\to m_{b}$. Moreover, by (I) we have $\ff_{(-\infty,0)}(v_{n})=\ff_{(-\infty,y_{n}^{-})}(w_{n})\leq\Lambda_{0}$ and $(i)$ follows by Lemma \ref{L:closel2}.\par\noindent
  Since by (\ref{eq:y+-}) we have $y_{n}^{+}-y_{n}^{-}\leq L_{0}$, by (II) we have $\ff_{(L_{0},+\infty)}(v_{n})= \ff_{(L_{0}+y_{n}^{-},+\infty)}(w_{n})\leq
 \ff_{(y_{n}^{+},+\infty)}(w_{n})\leq\Lambda_{0}$ and $(ii)$ follows by Lemma \ref{L:closel2}.\par\noindent
To prove $(iii)$ we first note that by Remark \ref{R:V-bound} we have $\|u\|^{2}\leq \tfrac{2\mu}{\mu-2}\beta$ for any
 $u\in \VV^{\beta}_{-}$.  Then, by $(i)$ we recover that $\|v_{n}(\cdot,y)\|^{2}_{2}\leq \tfrac{2\mu}{\mu-2}\beta+r_{0}$ for any $y\leq 0$.  Since $v_{n}(\cdot,0)=w_{n}(\cdot,y_{n}^{-})\in\VV^{\beta}_{-}$, by (\ref{eq:y+-}) we obtain moreover
 $\|v_{n}(\cdot,y)\|^{2}_{2}\leq \|v_{n}(\cdot,0)\|_{2}+\|v_{n}(\cdot,y)-v_{n}(\cdot,0)\|_{2}\leq\tfrac{2\mu}{\mu-2}\beta+ \tfrac{m_{b}+1}{\scriptstyle{\sqrt{2(b^{*}-b)}}}$ for any $y\in (0,y_{n}^{+}-y_{n}^{-}]$.
 Finally, by (II), we have $\|v_{n}(\cdot,y)-v_{n}(\cdot,y_{n}^{+}-y_{n}^{-})\|_{2}\leq 1$ for any $y>y_{n}^{+}-y_{n}^{-}$ and $(iii)$ follows with $\bar C=\tfrac{2\mu}{\mu-2}\beta+r_{0}+ \tfrac{m_{b}+1}{\scriptstyle{\sqrt{2(b^{*}-b)}}}+1$.\par\noindent
 To prove $(iv)$ we use $(iii)$ and (\ref{eq:Vmaggiore}). By (\ref{eq:Vmaggiore}) we know that there exists $C>0$ such that
\[\f(v_{n}(\cdot,y))\geq \tfrac 12\|\nabla v_{n}(\cdot,y)\|_{2}^{2} \left(1-C\tfrac{\| v_{n}(\cdot,y)\|_{2}^{(p+1)\theta}}{\|\nabla v_{n}(\cdot,y)\|_{2}^{2-(p+1)(1-\theta)}}\right)+
\tfrac 14\|v_{n}(\cdot,y)\|_{2}^{2}\quad\fa y\in\R.\] 
We set 
$$
{\cal A}_{n}=\{y\in\R\,|\,
\|\nabla v_{n}(\cdot,y)\|_{2}^{2-(p+1)(1-\theta)}\geq 2C\| v_{n}(\cdot,y)\|_{2}^{(p+1)\theta}\}.
$$
By (\ref{eq:Vmaggiore}), $\f(v_{n}(\cdot,y))\geq\tfrac 14\|v_{n}(\cdot,y)\|^{2}$ for every $y\in {\cal A}_{n}$ while $\|\nabla v_{n}(\cdot,y)\|_{2}^{2-(p+1)(1-\theta)}<2C\| v_{n}(\cdot,y)\|_{2}^{(p+1)\theta}$ for any $y\in\R\setminus {\cal A}_{n}$. By $(iii)$ we know that $\| v_{n}(\cdot,y)\|_{2}\leq \bar C$ for all $y\in\R$ and so  $\|\nabla v_{n}(\cdot,y)\|_{2}^{2}<\tilde C:= 2C\bar C^{(p+1)\theta}$ for any $y\in\R\setminus {\cal A}_{n}$. Given $(y_{1},y_{2})\subset\R$ we have
\begin{align*}
\|v_{n}\|_{H^{1}(S_{((y_{1},y_{2})})}^{2}&=\int_{y_{1}}^{y_{2}}\|\partial_{y}v_{n}(\cdot,y)\|_{2}^{2}+\|\nabla v_{n}(\cdot,y)\|^{2}_{2}+\|v_{n}(\cdot,y)\|_{2}^{2}\, dy\\
&\leq 2\ff(v_{n})+\int_{y_{1}}^{y_{2}}\|\nabla v_{n}(\cdot,y)\|^{2}_{2}\, dy+\bar C(y_{2}-y_{1})\\
&\leq 2\ff(v_{n})+\int_{(y_{1},y_{2})\cap {\cal A}_{n}}\|\nabla v_{n}(\cdot,y)\|^{2}_{2}\, dy+(\bar C+\tilde C)(y_{2}-y_{1})\\
&\leq 2\ff(v_{n})+4\int_{(y_{1},y_{2})\cap {\cal A}_{n}}\f(v_{n}(\cdot,y))-b\, dy+(\bar C+\tilde C+4b)(y_{2}-y_{1})\\
&\leq 6\ff(v_{n})+(\bar C+\tilde C+4b)(y_{2}-y_{1})\\&\leq \hat C^{2}=6(m_{0}+1)+(\bar C+\tilde C+4c)(y_{2}-y_{1})
\end{align*}
and $(iv)$ follows.\QED

By $(iv)$ of Lemma \ref{L:minseq} we have that the minimizing sequence $(v_{n})$ weakly converges in $H^{1}(S_{L})$ for any $L>0$ to a function $\bar v\in\MM$. Even if we do not know a priori that $\bar v\in\XX_{b}$, thanks to Lemma \ref{L:semicontinuitafun}, Lemma \ref{L:chisuradebole} and the semicontinuity of the distance function, the function $\bar v$ enjoyes the following properties  

\begin{Corollary}\label{C:minimo} For any $b\in [0,c)$ there exists $\bar v\in\MM$ such that
\begin{description}
\item{i)} given any interval $I\subset\R$ such that $\f(\bar v(\cdot,y))\geq b$ for a.e. $y\in I$ we have
$\ff_{I}(\bar v)\leq m_{b}$,
\item{ii)} $\mathrm{dist}(\bar v(\cdot,y),\VV^{\beta}_{-})\leq r_{0}$ for any $y\leq 0$,
\item{iii)} $\mathrm{dist}(\bar v(\cdot,y),\VV^{\beta}_{+})\leq r_{0}$ for any $y\geq L_{0}$,
\item{iv)} $\|\bar v(\cdot,y)\|_{2}\leq \bar C$ for any $y\in\R$,
\item{v)} for every $(y_{1},y_{2})\subset\R$, $\|\bar v\|_{H^{1}(S_{(y_{1},y_{2})})}\leq \hat C$,
\end{description}
where $L_{0}$, $\bar C$ and $\hat C$ are given by Lemma \ref{L:minseq}.
\end{Corollary}

We define $\bs$ and $\bt$ as follows:
\begin{align*}\bs&=\sup\{y\in \R\, /\, \dist(\bar v(\cdot,y),\VV^{b}_{-})\leq r_{0}\hbox{ and }
			\f(\bar v(\cdot,y))\leq
			b\},\\ \bt&=\inf\{y>\bs\, /\, 
			\f(\bar v(\cdot,y))\leq
			b\},
\end{align*}
with the agreement  that $\bs=-\infty$ whenever $\f(\bar v(\cdot,y))>b$
for every $y\in\R$ such that $\dist(\bar v(\cdot,y),\VV^{b}_{-})\leq r_{0}$ and that
$\bt=+\infty$ whenever $\f(\bar v(\cdot,y))>b$
for every $y>\bs$.

\begin{Remark}\label{R:tpm}{\sl Properties of $\bs$, $\bt$:} {\rm 
\begin{itemize}
\item[$i)$] $\bs\in [-\infty, L_{0}]$ and $\bt\in [0,+\infty]$.\par\noindent
By Corollary \ref{C:minimo}-($iii$), if $y\geq L_{0}$ then $\dist(\bar v(\cdot,y),\VV^{\ub}_{+})\leq r_{0}$. Hence $\dist(\bar v(\cdot,y),\VV^{b}_{-})> 4r_{0}$ for $y\geq L_{0}$ and $\bs\leq L_{0}$ follows. Moreover, by Corollary \ref{C:minimo}-($ii$), there results
$\dist(\bar v(\cdot,y),\VV^{\ub}_{-})\leq r_{0}$ if $y\leq 0$. Then, by the definition of $\bs$, we have that if $\bs<0$ then $\f(\bar v(\cdot,y))>b$ for any $y\in (\bs, 0]$ and so $\bt\geq 0$ follows. 
\item[$ii)$] If $\bs\in\R$ then $\bar v(\cdot,\bs)\in\VV^{b}_{-}$.\par\noindent
 Indeed, by definition, there exists a sequence $y_{n}\in(-\infty, \bs]$ such that $y_{n}\to \bs$ as $n\to+\infty$, $\f(\bar v(\cdot,y_{n}))\leq b$ and $\dist(\bar v(\cdot,y_{n}),\VV^{b}_{-})\leq r_{0}$ for any $n\in\N$. Then $\bar v(\cdot,y_{n})\in\VV^{b}_{-}$ for any $n\in\N$ and since, $\bar v(\cdot,y_{n})\to \bar v(\cdot,\bs)$ in $L^{2}(\R^{N})$, by Remark \ref{R:V-bound} we conclude that $\bar v(\cdot,\bs)\in\VV^{b}_{-}$.
\item[$iii)$] $\bs<\bt$. \par\noindent
It is sufficient to prove that if $\bs\in\R$ then, there exists $\delta>0$ such that $\f(\bar v(\cdot,y))>b$ for any $y\in (\bs,\bs+\delta)$. 
Assume by contradiction that there exists a sequence $(y_{n})\subset (\bs,+\infty)$ such that $\f(\bar v(\cdot,y_{n}))\leq b$ for any $n\in\N$ and $y_{n}\to \bs$. Then, by definition of $\bs$ we have $\dist(\bar v(\cdot,y_{n}),\VV^{b}_{-})> r_{0}$ for any $n\in\N$ and so 
$\bar v(\cdot,y_{n})\in\VV^{b}_{+}$. Hence, since $v(\cdot,y_{n})\to v(\cdot,\bar\s)$ in $L^{2}$, by Remark \ref{R:V-bound+}, we obtain $\bar v(\cdot,\bs)\in\VV^{b}_{+}$ while, by $(ii)$ we know that $\bar v(\cdot,\bs)\in\VV^{b}_{-}$.

\item[$iv)$] If $\bt\in\R$ then $\bar v(\cdot,\bt)\in\VV^{b}_{+}$.\par\noindent
 Indeed, by definition, there exists a sequence $y_{n}\in[\bt,+\infty)$ such that $y_{n}\to \bt$ as $n\to+\infty$, $\f(\bar v(\cdot,y_{n}))\leq b$. By definition of $\bs$, since $y_{n}>\bs$, we have $\dist(\bar v(\cdot,y_{n}),\VV^{b}_{-})> r_{0}$ for any $n\in\N$. Then $\bar v(\cdot,y_{n})\in\VV^{b}_{+}$ for any $n\in\N$ and since, $\bar v(\cdot,y_{n})\to \bar v(\cdot,\bt)$ in $L^{2}(\R)$, we conclude by Remark \ref{R:V-bound+} that $\bar v(\cdot,\bt)\in\VV^{b}_{+}$.

\item[$v)$] If $[y_{1},y_{2}]\subset  (\bs,\bt)$ then $\inf_{y\in [y_{1},y_{2}]}\f(\bar v(\cdot,y))>b$. Moreover $\ff_{(\bs,\bt)}(\bar v)\leq m_{b}$. \par\noindent
It follows by the definition of $\bs$ and $\bt$ that $\f(\bar v(\cdot,y))>b$ for any $y\in  (\bs,\bt)$. Then, by Lemma \ref{L:semicVv} we have $\inf_{y\in [y_{1},y_{2}]}\f(\bar v(\cdot,y))=\min_{y\in [y_{1},y_{2}]}\f(\bar v(\cdot,y))>b$ whenever $[y_{1},y_{2}]\subset  (\bs,\bt)$. By  Corollary \ref{C:minimo}--$(i)$ we furthermore derive that $\ff_{(\bs,\bt)}(\bar v)\leq m_{b}$.

\item[$vi)$] If $\bs=-\infty$ then $\liminf_{y\to-\infty}\f(\bar v(\cdot,y))-b=\liminf_{y\to-\infty}\dist(\bar v(\cdot,y),\VV^{b}_{-})=0$. \par\noindent
By  Corollary \ref{C:minimo}--$(ii)$ we have $\dist(\bar v(\cdot,y),\VV^{\ub}_{-})\leq r_{0}$ for every $y\leq 0$. Since $\bs=-\infty$ and $\ff_{(-\infty,\bt)}(\bar v)\leq m_{b}$ we derive that there exists a sequence
$y_{n}\to -\infty$ such that $\f(\bar v(\cdot,y_{n}))\to b$, $\bar v(\cdot,y_{n})\in\VV_{-}^{\beta}$ and $\dist(\bar v(\cdot,y_{n}),\VV^{b}_{+})\geq 4r_{0}$.\\ If $b=0$, by Remark \ref{R:V-bound},  we obtain $\bar v(\cdot,y_{n})\to 0$ and $(vi)$ follows. If $b>0$, arguing as in the proof of Lemma \ref{L:minimonormalim}, for any $n\in\N$, since $V(\bar v(\cdot,y_{n}))>b$, there exists of a unique $s_{n}\in (0,1]$ such that $\f( s_{n }\bar v(\cdot,y_{n}))=b$, $s_{n}\bar v(\cdot,y_{n}))\in\VV_{-}^{b}$ with $1-s_{n}\leq (\f(\bar v(\cdot,y_{n}))-b)/\nu^{-}(b)\to 0$. Since by Remark \ref{R:V-bound} $\|\bar v(\cdot,y_{n})||_{2}$ is bounded,  $\dist(\bar v(\cdot,y_{n}),\VV^{b}_{-})\leq (1-s_{n})\|\bar v(\cdot,y_{n})\|_{2}\to 0$ and $(vi)$ follows.
\item[$vii)$] If $\bt=+\infty$ then $\liminf_{y\to+\infty}\f(\bar v(\cdot,y))-b=\liminf_{y\to+\infty}\dist(\bar v(\cdot,y),\VV^{b}_{+})=0$. \par\noindent
The proof is analogous of the one of $(vi)$.
\end{itemize}
}
\end{Remark}
Thank to properties $(ii)$--$(iv)$, $(vi)$--$(vii)$ we recognize that the function $\bar v$ satisfies the assumption of Lemma \ref{L:st} on the interval $(\bs,\bt)$ which allows us to derive the following properties of $\bar v$. 
\begin{Lemma}\label{L:propv} There result
\begin{description}
\item{i)} $\ff_{(\bs,\bt)}(\bar v)=m_{b}$ and $\liminf_{y\to\bt^{-}}\f(\bar v(\cdot,y))=\liminf_{y\to\bs^{+}}\f(\bar v(\cdot,y))=b$,
\item{ii)} $\bt\in\R$ for any $b\in [0,c)$ and $\bs\in\R$ for any $b\in (0,c)$,
\item{iii)}  for every $h\in C_{0}^\infty(\R^{N}\times (\bs,\bt))$, with
 $\mathrm{supp}\, h\subset\R^{N}\times [y_{1},y_{2}]\subset \R^{N}\times(\bs,\bt)$, there exists $\bar t>0$ such that
 \begin{equation}\label{eq:minimal}\ff_{(\bs,\bt)}(\bar v+th)\geq\ff_{(\bs,\bt)}(\bar v),\quad\fa\, t\in (0,\bar t).\end{equation}
 Then
 $\bar v\in C^{2}(\R^{N}\times (\bs,\bt))$ verifies $-\Delta
    u+u-f(u)=0$ on $\R^{N}\times (\bs,\bt)$ and
    for any $[y_{1},y_{2}]\subset (\bs,\bt)$ there results
    $\bar v\in
    H^{2}(\R^{N}\times (y_{1},y_{2}))$,
\item{iv)} $E_{y}(\bar v(\cdot,y))=\tfrac 12\|\partial_{y}\bar v(\cdot,y)\|_{2}^{2}-\f(\bar v(\cdot,y))=-b\, $ for every $y\in (\bs,\bt)$,
\item{v)} $\liminf_{y\to\bt^{-}}\|\partial_{y}\bar v(\cdot,y)\|_{2}=\liminf_{y\to\bs^{+}}\|\partial_{y}\bar v(\cdot,y)\|_{2}=0$.
\end{description}
\end{Lemma}
\Proof
$(i)$ By Lemma \ref{L:st} we already know that $\ff_{(\bs,\bt)}(\bar v)\geq m_{b}$ and by $(v)$ of Remark \ref{R:tpm} we conclude that $\ff_{(\bs,\bt)}(\bar v)= m_{b}$. Hence, using Lemma \ref{L:st} again, we conclude that  $\liminf_{y\to\bt^{-}}\f(\bar v(\cdot,y))=\liminf_{y\to\bs^{+}}\f(\bar v(\cdot,y))=b$.\medskip\par\noindent
$(ii)$ Assume by contradiction that $\bt=+\infty$. 
By $(vii)$ of Remark \ref{R:tpm} there exists  $y_{0}>L_{0}$ such that
$u_{0}:=\bar v(\cdot,y_{0})\in\VV^{\bi_{+}}_{+}\setminus\VV^{ b}_{+}$ and $\bar v(\cdot,\cdot+y_{0})\in\XX_{b,u_{0}}^{+}$. To obtain a contradiction we show that 
\begin{equation}\label{eq:ccont}\fa y\geq y_{0}\hbox{ such that }V(\bar v(\cdot, y))\leq \bi_{+}\hbox{ we  have }\ff_{( y,+\infty)}(\bar v)\leq C_{+}(b)(\f(\bar v(\cdot, y))-b)^{3/2}.\end{equation}
By (\ref{eq:ccont}), using Lemma \ref{L:finitetime},  we derive that
there exists $\bar y\in (y_{0},y_{0}+1)$ such that $V(\bar v(\cdot, \bar y))=b$ which contradicts that $\bt=+\infty$.\\
If (\ref{eq:ccont}) does not hold, by Lemma \ref{L:minimonormalim}, 
there exists $\tilde y\geq y_{0}$ with $\bar v(\cdot,\tilde y)\in \VV^{\bi_{+}}_{+}$ and $\ff_{( \tilde y,+\infty)}(\bar v)> \ff_{( \tilde y,+\infty)}(w^{+}_{\bar v(\cdot, \tilde y)})$. Then, defining
$$
\tilde v(\cdot, y)=\begin{cases}\bar v(\cdot,y)&y\leq \tilde y\\ w^{+}_{\bar v(\cdot,\tilde y)}(\cdot,\cdot-\tilde y)&y> \tilde y
\end{cases}
$$ 
we obtain $\ff_{(\bs,+\infty)}(\tilde v)<\ff_{(\bs,+\infty)}(\bar v)=m_{b}$. On the other hand,  defining $\tilde \tau=\sup\{y>\bs\,|\,\f(\tilde v(\cdot,y))>b\}$, we recognize that $\tilde \tau$ satisfies the assumption of Lemma \ref{L:st} on the interval $(\bs,\tilde \tau)$ and we get the contradiction
$m_{b}\leq\ff_{(\bs,\tilde \tau)}(\tilde v)\leq \ff_{(\bs,+\infty)}(\tilde v)<m_{b}$. \\
To prove that $\bs\in\R$ when $b>0$ we can argue analogously using Lemmas \ref{L:minimonormalim-} and \ref{L:finitetime-}.\medskip\par\noindent
$(iii)$ 
Let us consider $h\in C_{0}^\infty(\R^{N}\times (\bs,\bt))$ with
 $\mathrm{supp}\, h\subset\R^{N}\times [y_{1},y_{2}]\subset \R^{N}\times(\bs,\bt)$.
By (v) of Remark \ref{R:tpm} we know that there exists $\mu>0$ such that $\f(\bar v(\cdot,y))\geq b+\mu$ for any $y\in [y_{1},y_{2}]$. Let us consider $(\bar v+th)^{*}$  the symmetric decreasing rearrangement of the function $v+th$ with respect to the variable $x$, i.e. the unique function with radial symmetry with respect to the variable $x\in\R^{N}$ such that 
$$
|\{x\in\R^{N}\,|\,(\bar v+th)^{*}(\cdot,y)> r\}|=|\{x\in\R^{N}\,|\, |(\bar v+th)(\cdot,y)|> r\}|\hbox{ for every $r>0$ and a.e. $y\in\R$}
$$
and $(\bar v+th)^{*}(x_{1},y)\geq (\bar v+th)^{*}(x_{2},y)$ whenever $|x_{1}|\leq |x_{2}|$,  for a.e. $y\in\R$.
 One recognizes (use e.g. \cite{[L]}, (12)-(14), and \cite{[LL]}, (3) pg. 73 ) that
$\|\nabla(\bar v+th)^{*}\|_{L^{2}(\R^{N}\times(y_{1},y_{2}))}\leq \|\nabla(\bar v+th)\|_{L^{2}(\R^{N}\times(y_{1},y_{2}))}$ and $\int_{\R^{N}\times(y_{1},y_{2})}\tfrac 12|(\bar v+th)^{*}|^{2}+F((\bar v+th)^{*})\, dx dy=
\int_{\R^{N}\times(y_{1},y_{2})}\tfrac 12|v+th|^{2}+F(|v+th|)\, dx dy=\int_{\R^{N}\times(y_{1},y_{2})}\tfrac 12|v+th|^{2}+F(\bar v+th)\, dx dy$. Therefore we have
\begin{align}
\nonumber\int_{\R^{N}\times [y_{1},y_{2}]}&\tfrac 12|\nabla(\bar v+th)^{*}|^{2}+
\tfrac 12|(\bar v+th)^{*}|^{2}-F((\bar v+th)^{*})\, dxdy\\ \label{eq:uno}
&\leq \int_{\R^{N}\times [y_{1},y_{2}]}\tfrac 12|\nabla(\bar v+th)|^{2}+
\tfrac 12|\bar v+th|^{2}-F(\bar v+th)\, dxdy\end{align}
We now claim that 
\begin{equation}\label{eq:mutilde}\exists\, \bar t>0\hbox{ such that }\f((\bar v+th)^{*}(\cdot,y))>b+\mu/2\hbox{ for any }t\in [0,\bar t]\hbox{ and }y\in [y_{1},y_{2}].\end{equation}
Arguing by contradiction, if (\ref{eq:mutilde}) does not hold, there exists a sequence $t_{n}\in (0,1)$ and a sequence $y_{n}\in[y_{1},y_{2}]$ such that $t_{n}\to 0$, $y_{n}\to y_{0}\in[y_{1},y_{2}]$ and $\f((\bar v+t_{n}h)^{*}(\cdot,y_{n}))\leq b+\mu/2$.\\ 
By Corollary \ref{C:minimo}-$(iv)$, since $h$ has compact support, We have that there exists $C>0$ such that 
$\| (\bar v+t_{n}h)^{*}(\cdot,y_{n})\|_{2}=\| (\bar v+t_{n}h)(\cdot,y_{n})\|_{2}\leq \|\bar v(\cdot,y_{n})\|_{2}+\|h(\cdot,y_{n})\|_{2}\leq C$ for any $n\in\N$. Since $\f((\bar v+t_{n}h)^{*}(\cdot,y_{n}))\leq b+\mu/2$, by Lemma \ref{L:coercivita}  there exists a constant $R>0$ such that $\|\nabla (\bar v+t_{n}h)^{*}(\cdot,y_{n})\|_{2}\leq R$ for any $n\in\N$. Then the sequence $\{(\bar v+t_{n}h)^{*}(\cdot,y_{n})\}$ is bounded in $H^{1}(\R^{N})$. Since the rearrangement is contractive in $L^{2}(\R^{N})$ we have  $\|(\bar v+t_{n}h)^{*}(\cdot,y_{n})-\bar v(\cdot,y_{0})\|_{2}\leq\|(\bar v+t_{n}h)(\cdot,y_{n})-\bar v(\cdot,y_{0})\|_{2}\to 0$ and so $(\bar v+t_{n}h)^{*}(\cdot,y_{n})\to \bar v(\cdot,y_{0})$ weakly in $H^{1}(\R^{N})$. By Lemma \ref{L:semiV} we then obtain the contradiction
$b+\mu/2\geq\liminf_{n\to+\infty}\f((\bar v+t_{n}h)^{*}(\cdot,y_{n}))\geq\f(\bar v(\cdot,y_{0}))\geq b+\mu$ which proves (\ref{eq:mutilde}).\\
Since $\bar v(\cdot,y)\in \X$ for a.e. $y\in\R$ we have $\bar v=\bar v^{*}$ and $(\bar v+th)^{*}=v$ for $x\in\R^{N}$ and $y\in\R\setminus [y_{1},y_{2}]$.
By (\ref{eq:mutilde}) we then recognize that  $(\bar v+th)^{*}$ satisfies the assumptions of Lemma \ref{L:st} on the interval $(\bs,\bt)$ for any $t\in [0,\bar t]$. Then $\ff_{(\bs,\bt)}((\bar v+t h)^{*})\geq m_{b}=\ff_{(\bs,\bt)}(\bar v)$
 for any $t\in [0,\bar t]$ and (\ref{eq:minimal}) follows by (\ref{eq:uno}). Finally, by (\ref{eq:minimal}) we have
 $$
 \int_{\R^{N}\times (\bs,\bt)} \tfrac 12|\nabla(\bar v+th)|^{2}+\tfrac 12|\bar v+th|^{2}-F(\bar v+th)-
 \tfrac 12|\nabla\bar v|^{2}-\tfrac 12|\bar v|^{2}+F(\bar v)\,dx dy\geq 0\quad \fa t\in (0,\bar t).
 $$
 Since $h$ is arbitrary  we derive that
 $\int_{\R^{N}\times (\bs,\bt)} \nabla\bar v\nabla h+\bar v\cdot h-f(\bar v)h\, dx\, dy=0$ for every $h\in C_{0}^{\infty}(\R^{N}\times (\bs,\bt))$,  and so that
 $\bar v$ is a weak solution of $(E)$ on  $\R^{N}\times (\bs,\bt)$.  Then $(iii)$ follows by $(v)$ of Corollary \ref{C:minimo} and standard regularity arguments.\medskip\par\noindent 
 $(iv)$ Fixed $\xi\in (\bs,\bt)$ and $s>0$ we define
$$\bar v_{s}(\cdot,y)=\begin{cases}\bar v(\cdot,y+\xi)&y\leq 0,\\
\bar v(\cdot,\tfrac{y}{s}+\xi)&0<y.
\end{cases}$$
and we note that $\bar v_{s}$ verifies the assumption of Lemma \ref{L:st}
on the interval $(\bs-\xi, s(\bt-\xi))$. Then 
$$
\ff_{(\bs-\xi, s(\bt-\xi))}(\bar v_{s})\geq m_{p}=\ff_{(\bs-\xi, \bt-\xi)}(\bar v(\cdot,\cdot+\xi))
$$
and so we have that for any $s>0$ there results
\begin{align*}
0&\leq\ff_{(\bs-\xi, s(\bt-\xi))}(\bar v_{s})-\ff_{(\bs-\xi, \bt-\xi)}(\bar v(\cdot,\cdot+\xi))\\
&=\int_{0}^{s(\bt-\xi)}\tfrac{1}{2}\|\partial_{y}\bar v_{s}(\cdot,y)\|^{2}+(\f(\bar v_{s}(\cdot,y))-b)\, dy-\int_{\xi}^{\bt}\tfrac{1}{2}\|\partial_{y}\bar v(\cdot,y)\|^{2}+(\f(\bar v(\cdot,y))-b)\, dy\\
&=\int_{0}^{s(\bt-\xi)}\tfrac{1}{2s^{2}}\|\partial_{y}\bar v(\cdot,\tfrac{y}{s}+\xi)\|^{2}+(\f(\bar v(\cdot,\tfrac{y}{s}+\xi))-b)\, dy-\ff_{(\xi, \bt)}(u)\\
&=\tfrac{1}{s}\int_{\xi}^{\bt}\tfrac{1}{2}\|\partial_{y}\bar v(\cdot,y)\|^{2}\, dy+
s\int_{\xi}^{\bt}\f(\bar v(\cdot, y))-b\, dy-\ff_{(\xi, \bt)}(\bar v)\\
&=(\tfrac{1}{s}-1)\int_{\xi}^{\bt}\tfrac{1}{2}\|\partial_{y}\bar v(\cdot,y)\|^{2}\, dy
+(s-1)\int_{\xi}^{\bt}\f(\bar v(\cdot, y))-b\, dy.
\end{align*}
This means that, setting $A=\int_{\xi}^{\bt}\tfrac{1}{2}\|\partial_{y}\bar v(\cdot,y)\|^{2}\, dy$ and $B=\int_{\xi}^{\bt}\f(\bar v(\cdot, y))-b\, dy$, the real function
$s\mapsto \psi(s)=A(\tfrac{1}{s}-1)+B(s-1)$ is non negative on $(0,+\infty)$ and then that
$
0\le \min \psi(s)=\psi({\scriptstyle\sqrt{\tfrac{A}{B}}})=
-(\sqrt{A}-\sqrt{B})^{2},
$
that implies $A=B$, i.e.,
\begin{equation}\label{eq:en}
\int_{\xi}^{\bt}\f(\bar v(\cdot,y))-b\, dy=\int_{\xi}^{\bt}\tfrac 12\|\partial_{y}\bar v(\cdot,y)\|_{2}^{2}\, dy\hbox{ for any }\xi\in (\bs,\bt).\end{equation} 
Since, by $(iii)$, $\bar v\in H^{2}(\R^{N}\times(y_{1},y_{2}))$ whenever $[y_{1},y_{2}]\subset (\bs,\bt)$, we derive that the function $y\to \tfrac 12\|\partial_{y}\bar v(\cdot,y)\|_{2}^{2}-V(\bar v(\cdot,y))$ is continuous and $(iv)$  follows by (\ref{eq:en}).
\medskip\par\noindent
$(v)$ It follows by $(i)$ and $(iv)$.\QED

\subsection{The case $b>0$. The periodic solutions} 
Consider the case $b\in (0,c)$. By (ii) of Lemma \ref{L:propv} we have $\bs, \bt\in\R$. In this case, by reflection and periodic conti\-nua\-tion,
starting from  $\bar v$, we can construct a solution to (E) on all $\R^{N+1}$ periodic in the variable $y$. Precisely let
$$
 v(x,y)=\begin{cases}\bar v(x,y+\bs)&\hbox{if }x\in\R^{N}\hbox{ and }y\in [0, \bt-\bs)\\
\bar v(x,\bt+(\bt-\bs-y))&\hbox{if }x\in\R^{N}\hbox{ and }y\in [\bt-\bs, 2(\bt-\bs)]\end{cases}
$$
and 
$v(x,y)= v(x,y+2k(\bt-\bs))$ for all $(x,y)\in\R^{N+1}$, $k\in\Z$.

\begin{Remark}\label{R:v}{\rm Let $T=\bt-\bs$. 
\begin{itemize}
\item[$i)$] The function $y\in\R\mapsto v(\cdot,y)\in L^{2}(\R^{N})$ is continuous and periodic with period $2T$. Moreover by $(ii)$ and $(iv)$ of Remark \ref{R:tpm},
$v(\cdot,0)\in\VV^{b}_{-}\hbox{ and }  v(\cdot,T)\in\VV^{b}_{+}.
$ Finally, by definition, $v(\cdot,-y)=v(\cdot,y)$ and $v(\cdot,y+T)=v(\cdot,T-y)$ for any $y\in\R$.

\item[$ii)$]  $v\in\MM$ and, by $(v)$ of Remark \ref{R:tpm}, $\f(v(\cdot,y))>b$ for any $y\in\R\setminus \{kT\,/\, k\in\Z\}$. 
\item[$iii)$] By $(v)$ of Lemma \ref{L:propv}, for any $k\in\Z$ we have
$\liminf_{y\to kT^{\pm}}\|\partial_{y}v(\cdot,y)\|_{2}=0.$
\item[$iv)$] By (iii) of Lemma \ref{L:propv} , $v\in C^{2}(\R^{N}\times (0, T))$ satisfies $-\Delta v(x,y)+v(x,y)-f(v(x,y))=0$ for $(x,y)\in \R^{N}\times (0, T)$.
\end{itemize}
}
\end{Remark}

We have

\begin{Lemma} \label{L:fin}
$v\in {\mathcal C}^{2}(\R^{N+1})$ is a solution of (E) on $\R^{N+1}$. Moreover,
$E_{v}(y)=\tfrac{1}{2}\|\partial_{y} v(\cdot,y)\|_{2}^{2}-\f(v(\cdot,y))=-b$ for all $y\in\R$ and
$\partial_{y} v(\cdot,0)=\partial_{y} v(\cdot,T)=0$. Finally $v>0$ on $\R^{N+1}$.
\end{Lemma}
\Proof First, let us prove that $v$ is a classical solution to 
(E). To this aim, we first note that by Remark  \ref{R:v} $(iii)$,
there exist four sequences $(\e^{\pm}_{n}), (\eta^{\pm}_{n})$, such that
$\e^{-}_{n}<0<\e^{+}_{n}$, $\eta^{-}_{n}<0<\eta^{+}_{n}$ for any $n\in\N$, $\e^{\pm}_{n},\, \eta^{\pm}_{n}\to 0$ and 
\begin{equation}\label{eq:sceltasequenze}\lim_{n\to+\infty}\|\partial_{y}v(\cdot,\e^{\pm}_{n})\|_{2}=\lim_{n\to+\infty}\|\partial_{y}v
(\cdot,T+\eta^{\pm}_{n})\|_{2}= 0.\end{equation}
Fixed any $\psi\in C_{0}^\infty(\R^{N+1})$, by Remark \ref{R:v} $(i)$-$(iv)$ we obtain that for any 
$k\in\Z$ and $n$ sufficiently large we have
\begin{align*}
    0& = 
    \int_{\R^{N}}\int_{2kT+\e_{n}^{+}}^{(2k+1)T+\eta_{n}^{-}}-\Delta v\,\psi+v\psi-f(v)\psi\,
    dy\, dx\\
    &=
     \int_{\R^{N}}\int_{2kT+\e_{n}^{+}}^{(2k+1)T+\eta_{n}^{-}}\nabla v\nabla\psi +v\psi-f(v)\psi\,
    dy\, dx +
    \int_{\R^{N}}\partial_{y}v(x,2kT+\e_{n}^{+})\psi(x,2kT+\e_{n}^{+})\, dx\\
    &\qquad-\int_{\R^{N}}\partial_{y}v(x,(2k+1)T+\eta_{n}^{-})\psi(x,(2k+1)T+\eta_{n}^{-})\, dx
        \end{align*}
        and
    \begin{align*}
    &0=
    \int_{\R^{N}}\int_{(2k-1)T+\eta_{n}^{+}}^{2kT+\e_{n}^{-}}-\Delta v\,\psi++v\psi-f(v)\psi\,
    dy\, dx\\
    &=
    \int_{\R^{N}}\int_{(2k-1)T+\eta_{n}^{+}}^{2kT+\e_{n}^{-}}\nabla v\nabla\psi ++v\psi-f(v)\psi\,
    dy\, dx-\int_{\R^{N}}\partial_{y}v(x,2kT+\e_{n}^{-})\psi(x,2kT+\e_{n}^{-})\, dx\\
    &\qquad +\int_{\R^{N}}\partial_{y}v(x,(2k-1)T+\eta_{n}^{+})\psi(x,(2k-1)T+\eta_{n}^{+})\, dx.
\end{align*}
By (\ref{eq:sceltasequenze}), in the limit for $n\to+\infty$, we obtain that for any $k\in\Z$ we have
$$0=\int_{\R^{N}}\int_{(2k-1)T}^{2kT}\nabla v\nabla\psi +v\psi-f(v)\psi\,
    dy\, dx=\int_{\R^{N}}\int_{2kT}^{(2k+1)T}\nabla v\nabla\psi +v\psi-f(v)\psi\,
    dy\, dx.$$
Then, $v$ satisfies 
$$
\int_{\R^{N+1}}\nabla v\nabla \psi++v\psi-f(v)\psi
     \, dx\, dy=0,\qquad
     \fa \psi\in C_{0}^\infty(\R^{N+1})
$$
and so $v$ 
is a {\sl classical} solution to
(E) on $\R^{N+1}$ which is periodic of period $2T$ in the variable $y$. Since by $(v)$ of Corollary \ref{C:minimo} we have $\|\bar v(\cdot,y)\|_{H^{1}(S_{(0,T)})}\leq\hat C$ depending only on $T$, by definition of $v$ and using (E) we recover that $v\in H^{2}(\R^{N}\times(y_{1},y_{2}))$ 
 for any bounded interval $(y_{1},y_{2})\subset \R$ and
 $\|v\|_{H^{2}(S_{(y_{1},y_{2})})}\leq C
$ with $C$ depending only on $y_{2}-y_{1}$.
This implies in particular that
the functions
$y\in \R\to\partial_{y}v(\cdot,y)\in L^{2}(\R^{N})$ and
$y\in \R\to v(\cdot,y)\in H^{1}(\R^{N})$ are uniformly 
continuous. Then 
$\lim_{y\to 0^{+}}\f(v(\cdot,y))-b=\liminf_{y\to 0^{+}}\|\partial_{y}v(\cdot,y)\|_{2}=0$ 
and analogously
$\lim_{y\to T^{-}}\f(v(\cdot,y))-b=\lim_{y\to T^{-}}\|\partial_{y}v(\cdot,y)\|_{2}=0$.
By continuity we derive that $\partial_{y} v(\cdot,0)=\partial_{y} v(\cdot,T)=0$.
By $(v)$ of Lemma \ref{L:propv} and the definition of $v$ it then follows that
$\tfrac12\|\partial_{y}v(\cdot,y)\|^{2}-\f(v(\cdot,y))=-b$ for any $y\in\R$. \\ To complete the proof we have to show that $v>0$ on $\R^{N+1}$. We know that $v\not=0$ and since $v\in\MM$ we have $v\geq 0$ on $\R^{N+1}$. Since $v$ solves (E) we have $-\Delta v+v=f(v)\geq 0$ on $\R^{N+1}$ and $v>0$  on $\R^{N+1}$ follows from the strong maximum principle.\QED

\begin{Lemma}\label{L:segnoder}
We have $\partial_{y}v>0$ on $\R^{N}\times (0,T)$ and  $\partial_{x_{i}}v<0$ on $\{(x,y)\in\R^{N+1}\,|\,
x_{i}>0\}$ for every $i=1,...,N$.\end{Lemma}
\Proof To prove that $\partial_{y}v>0$ on $\R^{N}\times (0,T)$ we first note that since
$\partial_{y}v(\cdot,0)=\partial_{y}v(\cdot,T)=0$ then $\partial_{y}v\in H^{1}_{0}(\R^{N}\times (0,T))$ and solves the linear elliptic equation $-\Delta\partial_{y}v+\partial_{y}v-f'(v)\partial_{y}v=0$ on
$\R^{N}\times (0,T)$. Then $\partial_{y}v\in H^{1}_{0}(\R^{N}\times (0,T))\cap H^{2}(\R^{N}\times (0,T))$ is an eigenfunction of the linear selfadjoint operator
$\mathcal{L}_{v}: H^{1}_{0}(\R^{N}\times (0,T))\cap H^{2}(\R^{N}\times (0,T))\subset L^{2}(\R^{N}\times (0,T))\to L^{2}(\R^{N}\times (0,T))$ defined by 
$\mathcal{L}_{v}h=-\Delta h+h-f'(v)h$ corresponding to the eigenvalue $0$.\\
The minimality property of $v$ proved in Lemma \ref{L:propv}-$(iii)$ implies $\langle\mathcal{L}_{v}h,h\rangle_{2}\geq 0$ for any $h\in C_{0}^{\infty}(\R^{N}\times (0,T))$ and we deduce that $0$ is the minimal eigenvalue of $\mathcal{L}_{v}$. Then $\partial_{y}v$ has constant sign on $\R^{N}\times (0,T)$. Assume by contradiction that $\partial_{y}v<0$ on $\R^{N}\times (0,T)$. Since, by construction, $v$ is even with respect to $T$, that implies that $v(x,T)\leq v(x,y)$ for all $x\in\R^{N}$ and $0<y<2T$. We deduce that $\partial^{2}_{y,y}v(x,T)\geq 0$ for all $x\in\R^{N}$ and so, multiplying (E) by $v$ and recalling that $v>0$ on $\R^{N+1}$ we deduce
$-\Delta_{x}v(x,T)\, v(x,T)+v(x,T)^{2}-f(v)v\geq 0$. Integrating with respect to $x$ on $\R^{N}$ we obtain
$\f'(v(\cdot,T))v(\cdot,T)\geq 0$ contrary to the fact that $v(\cdot,T)\in \VV^{b}_{+}$. This shows that
$\partial_{y}v>0$ on $\R^{N}\times (\bs,\bt)$. To prove that $\partial_{x_{i}}v<0$ on $\{(x,y)\in\R^{N+1}\,|\,
x_{i}>0\}$ we note that since $v\in \MM\cap C^{2}(\R^{N+1})$ we have
$\partial_{|x|}v(x,y)\leq 0$ for all $y\in\R$ and $|x|\not=0$. Then $\partial_{x_{i}}v\leq 0$ on $\{(x,y)\in\R^{N+1}\,|\, x_{i}>0\}$. Since $\partial_{x_{i}}v\leq 0$ solves the linear elliptic equation $-\Delta\partial_{x_{i}}v+\partial_{x_{i}}v=f'(v)\partial_{x_{i}}v$ on $\{(x,y)\in\R^{N+1}\,|\, x_{i}>0\}$ we deduce $-\Delta\partial_{x_{i}}v+\partial_{x_{i}}v\leq (f'(v))_{+}\partial_{x_{i}}v\leq 0$ on $\{(x,y)\in\R^{N+1}\,|\, x_{i}>0\}$ and since $\partial_{x_{i}}v\not=0$, the strong maximum principle assures $\partial_{x_{i}}v<0$ on $\{(x,y)\in\R^{N+1}\,|\, x_{i}>0\}$.\QED

\subsection{The case $b=0$. The homoclinic type mountain pass solution.} 

In the case $b=0$ Lemma \ref{L:propv} establishes  that $\bt\in\R$ but does not give information about $\bs$. We prove here below that in fact $\bs=-\infty$.
\begin{Lemma}\label{L:sigma} If $b=0$ then $\bs=-\infty$.\end{Lemma}
\Proof Assume that $\bs\in\R$. Then, arguing as in the case $b>0$, by reflection and periodic continuation, we construct a solution $v\in C^{2}(\R^{N+1})$ of (E) which is $2(\bt-\bs)$-periodic in the variable $y$ 
with $v(\cdot,0)\in\VV^{0}_{-}$ and $\partial_{ y} v(\cdot,0)=0$. Since $\VV^{0}_{-}=\{0\}$ we have
$v(x,0)=0$ and $\partial_{y}v(x,0)=0$ for any $x\in\R^{N}$.  Defining $a(x,y)=1-f(v(x,y))/v(x,y)$ when $v(x,y)\not= 0$ and $a(x,y)=1-f'(0)=1$ when $v(x,y)=0$ we have that $a$ is continuous on $\R^{N+1}$ and
$v$ solves $-\Delta v+a(x,y)v=0$ on $\R^{N+1}$. Defining the 
function $\tilde v(\cdot,y)= v(\cdot,y)$ for $y\in (0,2(\bt-\bs))$, and $\tilde v(\cdot,y)=0$ for $y\leq 0$ or $y\geq 2(\bt-\bs)$, since $v(x,0)=\partial_{y}v(x,0)=v(x,2(\bt-\bs))=\partial_{y}v(x,2(\bt-\bs))=0$, we obtain that also $\tilde v$ satisfies
$-\Delta v+a(x,y)v=0$ on $\R^{N+1}$. But a local unique continuation theorem (see e.g. Theorem 5 in \cite{[Ni]}) and a continuation argument imply that $\tilde v=0$ on $\R^{N+1}$ while, by definition of $\bs$ and $\bt$, $\tilde v(\cdot,y)= v(\cdot,y)=\bar v(\cdot,y+\bs)\not=0$ for $y\in (0,\bt-\bs)$.
 \QED
By Lemma \ref{L:sigma} we can define the function
$$
 v(x,y)=\begin{cases}\bar v(x,y+\bt)&\hbox{if }x\in\R^{N}\hbox{ and }y\in (-\infty,0]\\
\bar v(x,\bt-y)&\hbox{if }x\in\R^{N}\hbox{ and }y\in [0, +\infty)\end{cases}
$$
and the argument of the proof of Lemma \ref{L:fin} shows that
$v$ is a classical solution to (E) in $\R^{N+1}$. 
\begin{Remark}\label{R:Rfin}{\rm Again by $(v)$ of Corollary \ref{C:minimo} and using (E) we recover that $v\in H^{2}(\R^{N}\times(y_{1},y_{2}))$ 
 for any bounded interval $(y_{1},y_{2})\subset \R$ and
 $\|v\|_{H^{2}(S_{(y_{1},y_{2})})}\leq C
$ with $C$ depending only on $y_{2}-y_{1}$. This implies in particular that
the functions
$y\in \R\to\partial_{y}v(\cdot,y)\in L^{2}(\R^{N})$ and
$y\in \R\to v(\cdot,y)\in H^{1}(\R^{N})$ are uniformly 
continuous and so
$\lim_{y\to -\infty}\f(v(\cdot,y))=\liminf_{y\to +\infty}\f(v(\cdot,y))=0,$ 
$\lim_{y\to 0^{-}}\|\partial_{y}v(\cdot,y)\|_{2}=\liminf_{y\to 0^{+}}\|\partial_{y}v(\cdot,y)\|_{2}=0$, and 
 $E_{v}(y)=\tfrac 12\|\partial_{y}v(\cdot,y)\|_{2}^{2}-V(v(\cdot,y))=0$ for any $y\in\R$. Note finally that $v$ is radially symmetric with respect to $x$, not increasing with respect to $|x|$, and, by construction, even in the variable $y$.}\end{Remark}  

 In the case $b=0$ the functional $\ff(u)=\int_{\R}\tfrac 12\|\partial_{y}u(\cdot,y)\|^{2}_{2}+\f(u(\cdot,u))\, dy$ can be written, by Remark \ref{R:dimension}, as  $\ff(u)=\tfrac 12\|u\|^{2}_{H^{1}(\R^{N+1})}-\int_{\R^{N+1}}F(u)dxdy=\f_{N+1}(u)$ for all $u\in H^{1}(\R^{N+1})$ and in particular, denoting $c_{N+1}$ the mountain pass level of $\ff$ in $H^{1}(\R^{N+1})$, Proposition \ref{P:u0positive} establishes that $\ff$ has a positive radially symmetric critical point $w$ at the level $c_{N+1}$.  We  have  
 
\begin{Lemma}\label{L:mpmin} 
$v\in H^{1}(\R^{N+1})$  is a critical point of $\ff$ on $H^{1}(\R^{N+1})$ with $\ff(v)=c_{N+1}$. Moreover
$v\in {\mathcal C}^{2}(\R^{N})$ is a positive solution of (E) on $\R^{N+1}$ such that $v(x,y)\to 0$ as $|(x,y)|\to+\infty$, and, up to translations, $v$ is radially symmetric about the origin and $\partial_{r}v<0$ for $r=|(x,y)|>0$.\end{Lemma}
\Proof By Remark \ref{R:Rfin} we have $\lim_{y\to -\infty}\f(v(\cdot,y))=0$ and so there exists $y_{0}\leq-L_{0}<0$ ($L_{0}$ as in Corollary \ref{C:minimo}) such that $\f(v(\cdot,y))\leq\beta$ for any $y\leq y_{0}$. Since by Corollary \ref{C:minimo}-(ii) we know that  $\dist( v(\cdot,y),\VV^{\beta}_{+})\geq 4r_{0}$ for $y\leq -L_{0}$, we recognize that $v(\cdot,y)\in\VV^{\beta}_{-}$ for any $y\leq y_{0}$.  Then $\f'(v(\cdot,y))v(\cdot,y)\geq 0$ and by (\ref{eq:superV})  we obtain that $\f(v(\cdot,y))\geq\tfrac{\mu-2}{2\mu}\|v(\cdot,y)\|^{2}$ for any $y\leq y_{0}$. Since $m_{0}=\ff_{(-\infty,0)}(\bar v)\geq\int_{(-\infty,y_{0})}V( v(\cdot,y))\, dy$, using Corollary \ref{C:minimo}-(iv),
we then obtain
$\|v\|_{H^{1}(\R^{N+1})}^{2}=2\int_{-\infty}^{0}\|v(\cdot,y)\|^{2}\,dy\leq 2\int_{-\infty}^{y_{0}}\tfrac{2\mu}{\mu-2}\f(v(\cdot,y))\, dy+2\bar C|y_{0}|\leq \tfrac{4\mu}{\mu-2}m_{0}+2\bar C|y_{0}|$, 
and $v\in H^{1}(\R^{N+1})$ follows. Since $v\in H^{1}(\R^{N+1})$ solves (E) on $\R^{N+1}$, we deduce $v(x,y)\to 0$ as $|(x,y)|\to+\infty$, it is a critical point of $\ff$ on $H^{1}(\R^{N+1})$ and, by Remark \ref{R:tucritici}, $\ff(v)=2m_{0}\geq c_{N+1}$.\\
We now show that $2m_{0}\leq c_{N+1}$ and the Lemma will follows as in the proof of Proposition \ref{P:u0positive}.\\
As recalled above, $\ff$ admits on $H^{1}(\R^{N+1})$ a positive, radially symmetric (in $\R^{N+1}$) critical point $w$ such that $\ff(w)=c_{N+1}$ and $\partial_{r}w<0$ on $\R^{N+1}\setminus\{0\}$ where $r=|(x,y)|$. In particular $w(x,y)$ is radially symmetric with respect to
$x$ and monotone decreasing with respect to $|x|$  for any $y\in\R$ and so $w\in\MM$. By Lemma \ref{L:energia} we know that the {\sl energy} function $E_{w}(y)=\tfrac 12\|\partial_{y}w(\cdot,y)\|_{2}^{2}-\f(w(\cdot,y))$ is constant on $\R$. Since $w_{0}$ solves (E) we have $w\in H^{2}(\R^{N+1})\cap C^{2}(\R^{N+1})$. Then  $\|w(\cdot,y)\|\to 0$ and $\|\partial_{y}w(\cdot,y)\|_{2}\to 0$ as $y\to\pm\infty$ and we deduce that $E_{w}(y)=0$, i.e., $\tfrac 12\|\partial_{y}w(\cdot,y)\|_{2}^{2}=\f(w(\cdot,y))$ for any $y\in\R$. Since $w$ is even with respect to $y$ we have $\partial_{y}w(\cdot,0)=0$ and then $\f(w(\cdot,0))=0$. Since  $w$ is radially symmetric we have $w(\cdot,0)\not=0$ and so $w(\cdot,0)\in\VV^{0}_{+}$. Finally, since $\partial_{r}w<0$ on $\R^{N+1}\setminus\{0\}$ we derive that $\partial_{y}w(0,y)>0$ for any $y\in (-\infty,0)$ and we conclude
$\f(w(\cdot,y))=\tfrac 12\|\partial_{y}w(\cdot,y)\|_{2}^{2}>0$ for any $y\in (-\infty,0)$.\\
The above results tell us that $w$ satisfies the assumption of Lemma \ref{L:st}
on the interval $(-\infty,0)$ and $\ff_{(-\infty,0)}(w)\geq m_{0}$ follows. Hence $c_{N+1}=\ff(w)\geq 2m_{0}$ and we conclude that $c_{N+1}=2m_{0}$.\\
To conclude the proof we note that since $v\geq 0$ and $-\Delta v+v=f(v)\geq 0$ on $\R^{N+1}$, the strong maximum principle establishes that
$v>0$ on $\R^{N+1}$ and so, by Theorem 1 in \cite{[LN]} we conclude that, up to translations, $v$ is radially symmetric  about the origin and $\partial_{r}v<0$ for $r=|(x,y)|>0$\QED

\end{document}